\documentclass[preprint]{elsarticle}
\usepackage{eurosym}
\usepackage[english]{babel}
\usepackage[utf8]{inputenc}
\usepackage[T1]{fontenc}
\usepackage{amsmath}
\usepackage{amsfonts}
\usepackage{amsthm}
\usepackage{amssymb}
\usepackage{graphicx}
\usepackage{dsfont}
\usepackage{color}

\newtheorem{theorem}{Theorem}

\newtheorem{lemma}[theorem]{Lemma}

\begin{document}

\begin{frontmatter}
\title{A Mixed Finite Element Method for a Class of Evolution Differential Equations with $p$-Laplacian and Memory}

\author[label1]{Rui M.P. Almeida}
\ead{ralmeida@ubi.pt}

\author[label1]{Jos\'{e} C.M. Duque}
\ead{jduque@ubi.pt}

\author[label1]{Belchior C.X. M\'{a}rio}
\ead{belchior.mario@ubi.pt}

\affiliation[label1]{organization={University of Beira Interior, Center of Mathematics and Applications},
             city={Covilh\~a},
             country={Portugal}}

\begin{abstract}
We present a new mixed finite element method for a class of parabolic equations with $p$-Laplacian and nonlinear memory. The applicability, stability and convergence of the method are studied. First, the problem is written in a mixed formulation as a system of one parabolic equation and a Volterra equation. Then, the system is discretized in the space variable using the finite element method with Lagrangian basis of degree $r\geq1$. Finally, the Cranck-Nicolson method with the trapezoidal quadrature is applied to discretize the time variable. For each method, we establish existence, uniqueness and regularity of the solutions. The convergence order is found to be dependent on the parameter $p$ on the $p$-Laplacian in the sense that it decreases as $p$ increases.
\end{abstract}


\end{frontmatter}

\section{Introduction}

In this work, we study the evolutionary integro-differential equation with $%
p $-Laplacian and memory,
\begin{equation}  \label{eq:30}
u_{t}(x,t)-\Delta_{p}u(x,t)=\int_{0}^tg(t-s)\Delta_{p}u(x,s)ds+f(x,t),
\end{equation}
where $f$ and $g$ are given functions. This type of equation appears in the
mathematical description of heat propagation in materials with memory, where
the heat flux can depend on the process history. In addition to the common
problem of the memory term, whose time discretization forces a large volume
of calculations and memory consumption, problem (\ref{eq:30}) also presents
a new difficulty when containing the $p$-Laplacian in that it makes the
memory term nonmonotonous (see \cite{MR3462616}).

Since the 70's, evolution equations with memory terms have attracted the
attention of researchers. Issues related to the existence and properties of
solutions to partial integro-differential equations (PIDEs) such as
\begin{equation}  \label{eqlm}
u_{t}-Au=\int_{0}^{t}g(t-s)Bu(s)ds+f,
\end{equation}%
where $A$ and $B$ are symmetric positive definite operators of at most
second order and $g$ is a memory kernel, were addressed, for example, in
\cite{BM79,CLN78,Cam76,Noh79}. Concerning numerical approximations to the
solutions of (\ref{eqlm}), several methods have already been investigated.
For a review of the finite element method applied to PIDEs, we refer, for
instance, to \cite{chen_book} and the references therein. The finite volume
method and the collocation method with splines were studied in \cite{SEL06}
and \cite{PFF08}, respectively. The mixed finite element method \ was
considered in \cite{SEL09}, while the discontinuous Galerking method was
studied in \cite{MBMS11}. Subsequently, an error analysis for the
Cranck-Nicolson finite element method was made in \cite{RSC19}. A two grid
finite element method was proposed in \cite{WH19} and, more recently, a
pseudospectral method was investigated in \cite{TDTR21}.

Partial differential equations involving the $p$-Laplacian operator have
been extensively studied in the last decades. For a survey of the theory, we
refer to the monographs \cite{AS_book,DiB_book,Lio_book}. Concerning the
numerical simulations of the $p$-Laplacian with the finite element method,
it was found that the regularity of the solutions limits the convergence
rates. In \cite{GM75}, Glowinski and Marroco proved a convergence of $%
\mathcal{O}(h^{\frac{1}{p-1}})$ in the $W^{1,p}$ norm. Later, Chow \cite%
{Cho89} improved this convergence order to $\mathcal{O}(h^{\frac{2}{p}})$.
In 1993, assuming a stronger regularity for the weak solution, Barret and
Liu, in \cite{BL94}, proved optimal error bounds of order $\mathcal{O}(h)$
in $W^{1,p}$.

The lack of monotonicity in the memory term of equation (\ref{eq:30}) makes
it unfeasible to use most of the well-developed techniques available. In
\cite{MR3462616}, Antontsev and his coauthors studied equation (\ref{eq:30})
with a nonlinear source term $\Theta (x,t,u)$, substituting the equation with
a system composed of a diffusion-reaction equation and an integral equation.
They proved that for $\max \{1,\frac{2n}{n+2}\}<p<\infty $, $u_{0}\in
W_{0}^{1,p}(\Omega ),f\in L^{2}(Q)$ and $g,g^{\prime}\in L^2(0,T)$, the problem
has a weak solution that is local or global in time depending on the growth rate of $%
\Theta (x,t,s)$, when $\left\vert s\right\vert \rightarrow \infty $.
Uniqueness conditions were established and they also proved that for $p>2$
and $s\Theta (x,t,s)\leq 0$, the data disturbances propagate with finite
speed and that the waiting time effect is possible.

Nowadays, problem (\ref{eq:30}), with $p$ depending on $x$, is attracting
considerable attention, perhaps because of its various physical
applications. We refer to \cite{MR3919666,ZM20,MR4220412}, where questions
on the solvability and properties of the solutions are addressed.

In this paper, we present a new mixed finite element method for equation (%
\ref{eq:30}). The existence, uniqueness and regularity of the discrete
solutions are established. Error bounds depending on the parameter $p$ are
also obtained. An auxiliary problem and its variational formulation is
presented in Section 2. The discretization of the space variable is
developed in Section 3. The discretization of the time variable is studied
in Section 4. Finally, in Section 5, we draw some final conclusions.

\section{Parabolic equation with $p$-Laplacian}

Let us consider the evolutionary integro-differential equation with the
homogeneous Dirichlet condition,
\begin{eqnarray}  \label{eq:14}
\begin{cases}
u_{t}-\Delta_{p}u=\int_{0}^{t}g(t-s)\Delta_{p}u(x,s)ds+f(x,t),\,\forall(x,t)%
\in Q=\Omega\;\times]0,T], \\
u(x,t)=0,\,\forall(x,t)\in\partial\Omega\times[0,T], \\
u(x,0)=u_{0}(x),\,\forall x\in\Omega,%
\end{cases}%
\end{eqnarray}
where $u_0$, $g$ and $f$ are given functions, $\Omega\subset\mathbb{R}^n$ is
a bounded domain with Lipschitz-continuous boundary. The $p$-Laplacian $%
\Delta_{p}u$ is given by
\begin{equation*}
\Delta_{p}u=\mathrm{div}\left(\left|\nabla u\right|^{p-2}\nabla
u\right),\quad 2<p<\infty
\end{equation*}
and
\begin{equation}  \label{defy}
y(x,t)=\int_{0}^{t}g(t-s)\Delta_{p}u(x,s)ds
\end{equation}
is the memory term of the evolutionary integro-differential equation.

Assuming
\begin{equation*}
g,g^{\prime}\in L^2(0,T),\quad f\in L^{2}(Q),\quad u_{0}\in L^{2}(\Omega )\cap
W_{0}^{1,p}(\Omega ),
\end{equation*}%
it is proved in \cite{MR3462616} that problem (\ref{eq:14}) has a unique
weak solution. For notions on Sobolev spaces, we refer to \cite{MR1801735,
MR1625845,Lio_book}. Below, we present an important lemma which can be found
in \cite{BL94,MR3250760,Cho89}.

\begin{lemma}
\label{eq:le1} For every $p>1$ and $\delta\geq 0$, there are three positive
constants $C_1$, $C_2$ and $C_3$ such that for every $\zeta,\gamma\in\mathbb{%
R}^n$, $\zeta\neq\gamma$, we have

\begin{enumerate}
\item
\begin{equation*}
\left|\left|\zeta\right|^{p-2}\zeta-\left|\gamma\right|^{p-2}\gamma\right|%
\leq C_1\left|\zeta-\gamma
\right|^{1-\delta}\left(\left|\zeta\right|+\left|\gamma\right|\right)^{p-2+%
\delta};
\end{equation*}

\item
\begin{equation*}
\left(\left|\zeta\right|^{p-2}\zeta-\left|\gamma\right|^{p-2}\gamma,\zeta-%
\gamma\right)_{\mathbb{R}^n}\geq
C_2\left|\zeta-\gamma\right|^{2+\delta}\left(\left|\zeta\right|+\left|\gamma%
\right|\right)^{p-2-\delta};
\end{equation*}

\item if $p>2$,
\begin{equation*}
\left(\left|\zeta\right|^{p-2}\zeta-\left|\gamma\right|^{p-2}\gamma,\zeta-%
\gamma\right)_{\mathbb{R}^n}\geq C_3\left|\zeta-\gamma\right|^{p}.
\end{equation*}
\end{enumerate}
\end{lemma}

\subsection{Auxiliary problem}

It is easy to prove that the memory term satisfies the integral equation
\begin{eqnarray}  \label{eq:jaa}
y(x,t)=-\int_{0}^{t}g(t-s)y(x,s)ds+f_2(x,t,u),
\end{eqnarray}
where $f_2$ is the nonlinear nonlocal operator
\begin{eqnarray}  \label{eq:f24}
f_2(x,t,u(x,t))&=&u(x,t)g(0)-u_0(x)g(t)+\int_{0}^{t}g^{\prime }(t-s)u(x,s)ds
\notag \\
&&-\int_{0}^{t}g(t-s)f(x,s)ds.
\end{eqnarray}
In fact, taking equation (\ref{eq:14}) and convoluting with $g$, we obtain
\begin{eqnarray}  \label{eq:hh}
&&\int_{0}^{t}g(t-s)\Delta_pu(x,s)ds  \notag \\
&&=-\int_{0}^{t}g(t-s)\int_{0}^{s}g(s-\tau)\Delta_pu(x,\tau)d\tau
ds+\int_{0}^{t}g(t-s)u_s(x,s)ds  \notag \\
&&-\int_{0}^{t}g(t-s)f(x,s)ds.
\end{eqnarray}
Integrating by parts the second term on the right side of equation (\ref%
{eq:hh}), we obtain
\begin{eqnarray*}
\int_{0}^{t}g(t-s)\Delta_pu(x,s)ds&=&-\int_{0}^{t}g(t-s)\int_{0}^{s}g(s-%
\tau)\Delta_p u(x,\tau)d\tau ds  \notag \\
&&+g(0)u(x,t)-g(t)u_0(x)+\int_{0}^{t}g^{\prime }(t-s)u(x,s)ds  \notag \\
&&-\int_{0}^{t}g(t-s)f(x,s)ds,
\end{eqnarray*}
which means that $y(x,t)$, defined in (\ref{defy}), satisfies equation (\ref%
{eq:jaa}).

This allows us to consider the equivalent auxiliary problem of finding the
pair $(u,y)$ that satisfies the conditions
\begin{eqnarray}  \label{eq:ff23}
\begin{cases}
u_{t}-\Delta_{p}u=y(x,t)+f(x,t),\quad\forall(x,t)\in Q, \\
y(x,t)=-\int_{0}^{t}g(t-s)y(x,s)ds+f_2(x,t,u),\quad\forall(x,t)\in Q, \\
u(x,t)=0,\quad\forall(x,t)\in\partial\Omega\times[0,T], \\
u(x,0)=u_{0}(x),\quad\forall x\in\Omega, \\
y(x,0)=0,\quad\forall x\in\Omega.%
\end{cases}%
\end{eqnarray}

\subsection{Variational formulation}

If we multiply the first equation of problem (\ref{eq:ff23}) by $w\in
H_{0}^{1}(\Omega)$ and integrate in $\Omega$, we get
\begin{eqnarray}  \label{eq:2}
\int_{\Omega}^{}u_{t}wdx-\int_{\Omega}^{}\mathrm{div}\left(\left|\nabla
u\right|^{p-2}\nabla u\right)wdx =\int_{\Omega}^{}ywdx+\int_{\Omega}^{}fwdx.
\end{eqnarray}

Applying Green's theorem to the second term on the right side of
equation (\ref{eq:2}) and using the definition of the space $%
H_{0}^{1}(\Omega)$, we have
\begin{eqnarray*}  \label{eq:5}
\int_{\Omega}^{}u_{t}wdx+\int_{\Omega}^{}\left|\nabla u\right|^{p-2}\nabla
u\nabla wdx=\int_{\Omega}^{}ywdx+\int_{\Omega}^{}fwdx.
\end{eqnarray*}

Multiplying the second equation of the problem (%
\ref{eq:ff23}) by $v\in H_{0}^{1}(\Omega)$  and integrating in $\Omega$ we obtain
\begin{eqnarray*}  \label{eq:js}
\int_{\Omega}^{}yvdx=-\int_{\Omega}^{}v\int_{0}^{t}g(t-s)y(x,s)dsdx+\int_{%
\Omega}^{}f_2(x,t,u)vdx.
\end{eqnarray*}

So, we are left with the system
\begin{eqnarray}  \label{eq:ff27}
\begin{cases}
\int_{\Omega}^{}u_twdx+\int_{\Omega}^{}\left|\nabla u\right|^{p-2}\nabla
u\nabla wdx=\int_{\Omega}^{}ywdx+\int_{\Omega}^{}fwdx,\,\forall w\in
H_0^1(\Omega), \\
\int_{\Omega}^{}yvdx=-\int_{\Omega}^{}v\int_{0}^{t}g(t-s)y(x,s)dsdx+\int_{%
\Omega}^{}f_2vdx,\,\forall v\in H_0^1(\Omega).%
\end{cases}%
\end{eqnarray}

A pair of bounded and measurable functions $(u(x,t),y(x,t))$, is said to be a
weak solution of the initial value problem (\ref{eq:ff23}), with initial
data $u_0$ limited and measurable, if (\ref{eq:ff27}) is valid for all $%
(w,v)\in(H_0^1(\Omega))^2$. Henceforth, we assume that problem (\ref{eq:ff23}%
) has a unique weak solution with sufficient regularity in order to perform the
calculations needed in next sections.

\section{Discretization in space}

\subsection{Lagrangian Bases}

Let us consider a regular partition $\mathcal{T}_h=\{T_0,\dots,T_m\}$ of $%
\Omega$ in simplexes with parameter $h$ and the space $\mathcal{S}^h\subset
H_0^1(\Omega)$ defined by
\begin{equation*}
\mathcal{S}^h=\{w\in C^0(\Omega):w(x)=0,\;x\in\partial\Omega,\;w(x)|_{T_k}\in%
\mathcal{P}_r(T_k),\;k=0,\cdots,m\},
\end{equation*}
where $\mathcal{P}_r(T_k)$ is the set of polynomials of degree less than or
equal to $r$ defined in $T_k$. We denote the interpolation operator into $\mathcal{S}^h$ by $\Pi_h$.
An estimate of the interpolation error is given in the next Lemma, which may be found in \cite{MR1930132}.

\begin{lemma}
\label{eq:le3} If $\Pi_h:H^{r+1}(\Omega)\cap H_0^1(\Omega)\rightarrow%
\mathcal{S}^h$ is the interpolation operator, then
\begin{equation*}
\|u-\Pi_hu\|_{L^2(\Omega)}+h\|\nabla(u-\Pi_hu)\|_{L^2(\Omega)}\leq
Ch^s\|u\|_{H^s(\Omega)},\quad 1\leq s\leq r+1,
\end{equation*}
with $u\in H^s(\Omega)\cap H_0^1(\Omega)$ and $C$ is a positive constant.
\end{lemma}

The semi-discrete problem is to find $(u^h,y^h)\in \big(\mathcal{S}^h\big)^2$%
such that
\begin{eqnarray}  \label{eq:ff28}
\begin{cases}
\int_{\Omega}^{}u_t^hw^hdx+\int_{\Omega}^{}\left|\nabla
u^h\right|^{p-2}\nabla u^h\nabla
w^hdx=\int_{\Omega}^{}y^hw^hdx+\int_{\Omega}^{}fw^hdx,\,\forall w^h\in%
\mathcal{S}^h, \\
\int_{\Omega}^{}y^hv^hdx=-\int_{\Omega}^{}v^h\int_{0}^{t}g(t-s)y^h(x,s)dsdx+%
\int_{\Omega}^{}f_2^hv^hdx,\,\forall v^h\in\mathcal{S}^h,%
\end{cases}%
\end{eqnarray}
and
\begin{equation*}
u^h(x,0)=u_0^h=\Pi_hu_0,\;\,y^h(x,0)=0,\;\,\forall x\in\Omega.
\end{equation*}
We note that problem (\ref{eq:ff28}) has a solution. In fact, from the second equation we
obtain a solution $y^h(u^h)$ which is substituted in the first equation to  give a
solution $u^h(x,t)$. Substituting  the latter solution in $y^h(u^h)$ gives the
solution $y^h(x,t)$ (see \cite{MR3462616} for more details).

\begin{theorem}[Uniqueness]
\label{teo_uni} If $g,g^{\prime}\in L^{\infty}(0,T)$, then the solution of the
semi-discrete problem (\ref{eq:ff28}) is unique.
\end{theorem}

\begin{proof}
Suppose that $(u_1^h,y_1^h)$ and $(u_2^h,y_2^h)$ are two solutions of the semi-discrete problem (\ref{eq:ff28}). Subtracting the equation for $y_2^h$ from the equation for $y_1^h$, we obtain
\begin{eqnarray}
\label{eq:js1}
\int_{\Omega}^{}(y_1^h-y_2^h)v^hdx&=&-\int_{\Omega}^{}v^h\int_{0}^{t}g(t-s)(y_1^h-y_2^h)dsdx +g(0)\int_{\Omega}^{}(u_1^h-u_2^h)v^hdx\nonumber\\
&&+\int_{\Omega}^{}v^h\int_{0}^{t}g'(t-s)(u_1^h-u_2^h)dsdx.
\end{eqnarray}
Let $v^h=y_1^h-y_2^h\in \mathcal{S}^h$. As $g$ and $g'$ are bounded, we may apply Young's inequality to (\ref{eq:js1}) and thus obtain
\begin{eqnarray}
\label{eq:js2}
\int_{\Omega}^{}(y_1^h-y_2^h)^2dx&\leq &C\int_{0}^{t}\int_{\Omega}^{}(y_1^h-y_2^h)^2dxds+C\int_{\Omega}^{}(u_1^h-u_2^h)^2dx\nonumber\\&&+C\int_{0}^{t}\int_{\Omega}^{}(u_1^h-u_2^h)^2dxds.
\end{eqnarray}
Applying Gronwall's lemma to (\ref{eq:js2}), we have
\begin{eqnarray}
\label{eq:js3}
\int_{\Omega}^{}(y_1^h-y_2^h)^2dx\leq C\int_{\Omega}^{}(u_1^h-u_2^h)^2dx+C\int_{0}^{t}\int_{\Omega}^{}(u_1^h-u_2^h)^2dxds.
\end{eqnarray}
Using a similar argument with $u_1^h$ and $u_2^h$, we get
\begin{eqnarray}
\label{eq:js4}
&&\int_{\Omega}^{}\big((u_1^h)_t-(u_2^h)_t\big)w^hdx+\int_{\Omega}^{}\left(\left|\nabla u_1^h\right|^{p-2}\nabla u_1^h-\left|\nabla u_2^h\right|^{p-2}\nabla u_2^h\right)\nabla w^hdx\nonumber\\&&=\int_{\Omega}^{}(y_1^h-y_2^h)w^hdx.
\end{eqnarray}
Let $w^h=u_1^h-u_2^h\in \mathcal{S}^h$. Applying Young's inequality in (\ref{eq:js4}), using (\ref{eq:js3}) and Lemma \ref{eq:le1}, we have
\begin{eqnarray}
\label{eq:js5}
&&\frac{d}{dt}\int_{\Omega}^{}(u_1^h-u_2^h)^2dx\nonumber\\&&\leq\int_{\Omega}^{}(u_1^h-u_2^h)^2dx+C\int_{\Omega}^{}(u_1^h-u_2^h)^2dx+C\int_{0}^{t}\int_{\Omega}^{}(u_1^h-u_2^h)^2dxds.
\end{eqnarray}
Integrating equation (\ref{eq:js5}) from $0$ to $t$ and noting that $u_1^h(x,0)-u_2^h(x,0)=0$, we obtain
\begin{eqnarray*}
\label{eq:js6}
\int_{\Omega}^{}(u_1^h-u_2^h)^2dx\leq C(1+t)\int_{0}^{t}\int_{\Omega}^{}(u_1^h-u_2^h)^2dxds.
\end{eqnarray*}
Now, let $\zeta(t)=\int_{\Omega}^{}(u_1^h-u_2^h)^2dx\geq 0$. This implies that $\zeta(t)\leq C\int_{0}^ {t}\zeta(s)ds$ and $\zeta(0)=0$. By Gronwall's lemma, $\zeta(t)=0$ for every $t\in [0,T]$, that is,
\begin{eqnarray*}
\label{eq:js8}
\int_{\Omega}^{}(u_1^h-u_2^h)^2dx=0,\quad u_1^h=u_2^h\;\,\mathrm{in}\;\,L^2(\Omega).
\end{eqnarray*}
Returning to $y$, from equation (\ref{eq:js3}), we have
\begin{equation*}
\int_{\Omega}^{}(y_1^h-y_2^h)^2dx=0,\quad\mathrm{for\;every}\;t\in[0,T],
\end{equation*}
which proves the desired result.
\end{proof}

\begin{theorem}[Regularity]
Let $(u^h,y^h)\in\mathcal{S}^h\times\mathcal{S}^h$ be a solution to
problem (\ref{eq:ff28}) and $g,g^{\prime} \in L^{\infty}(0,T)$.  Then, for every $ t \in
[0,T]$,
\begin{equation}  \label{reg_u}
\|u^h\|_{L^2(\Omega)}^2+\|\nabla u^h\|_{L^p(Q)}^p\leq
C\|u_0\|_{L^2(\Omega)}^2+C\|f\|_{L^2(Q)}^2
\end{equation}
and
\begin{equation*}
\|y^h\|_{L^2(\Omega)}^2\leq C\|u_0\|_{L^2(\Omega)}^2+C\|f\|_{L^2(Q)}^2,
\end{equation*}
where
\begin{equation*}
C=C\big(T,\|g\|_{L^\infty(0,T)},\|g^{\prime }\|_{L^\infty(0,T)}\big).
\end{equation*}
\end{theorem}

\begin{proof}
Since $y^h$ is a solution to problem (\ref{eq:ff28}),
\begin{eqnarray*}
\label{eq:mm1}
\int_{\Omega}^{}y^hv^hdx=-\int_{\Omega}^{}v^h\int_{0}^{t}g(t-s)y^h(x,s)dsdx+\int_{\Omega}^{}f_2v^hdx.
\end{eqnarray*}
For $v^h=y^h\in\mathcal{S}^h$,
\begin{eqnarray}
\label{eq:mm2}
\int_{\Omega}^{}(y^h)^2dx&=&-\int_{\Omega}^{}y^h(x,t)\int_{0}^{t}g(t-s)y^h(x,s)dsdx+g(0)\int_{\Omega}^{}u^hy^hdx\nonumber\\
&&-g(t)\int_{\Omega}^{}u_0^hy^hdx+\int_{0}^{t}g'(t-s)\int_{\Omega}^{}u^h(x,s)y^h(x,t)dxds\nonumber\\
&&-\int_{0}^{t}g(t-s)\int_{\Omega}^{}f(x,s)y^h(x,t)dxds.
\end{eqnarray}
Applying Young's inequality to equation (\ref{eq:mm2}) and using the fact that $g$ and $g'$ are limited, we have
\begin{eqnarray}
\label{eq:mm8}
\int_{\Omega}^{}(y^h)^2dx&\leq& C\int_{0}^{t}\int_{\Omega}^{}\big(y^h(x,s)\big)^2dxds+C\int_{\Omega}^{}(u^h)^2dx +C\int_{\Omega}^{}(u^h_0)^2dx\nonumber\\
&&+\int_{0}^{t}\int_{\Omega}^{}(u^h)^2dxds+C\int_{0}^{t}\int_{\Omega}^{}f^2dxds.
\end{eqnarray}
Applying Gronwall's lemma to (\ref{eq:mm8}), we obtain
\begin{eqnarray}
\label{eq:set}
\|y^h\|_{L^2(\Omega)}^2 \leq C\|u^h\|_{L^2(\Omega)}^2+C\|u^h_0\|_{L^2(\Omega)}^2+C\int_{0}^{t}\|u^h\|_{L^2(\Omega)}^2ds+C\|f\|_{L^2\left(Q\right)}^2.&&
\end{eqnarray}
Now, considering the first equation of problem (\ref{eq:ff28}) with $w^h=u^h$, we have
\begin{eqnarray}
\label{eq:mm3}
\int_{\Omega}^{}u^h_tu^hdx+C\int_{\Omega}^{}|\nabla u^h|^{p}dx=\int_{\Omega}^{}y^hu^hdx+\int_{\Omega}^{}fu^hdx.
\end{eqnarray}
Applying Young's inequality to (\ref{eq:mm3}) and equation (\ref{eq:set}) we have
\begin{eqnarray}
\label{eq:mm5}
\frac{d}{dt}\|u^h\|^2_{L^2(\Omega)}+C\|\nabla u^h\|^p_{L^p(\Omega)}&\leq& C\|u^h\|^2_{L^2(\Omega)}+C\|u_0^h\|^2_{L^2(\Omega)} +C\int_{0}^{t}\|u^h\|^2_{L^2(\Omega)}ds\nonumber\\ &&+C\|f\|^2_{L^2\big(0,T;L^2(\Omega)\big)}+\frac{1}{2}\|f\|^2_{L^2(\Omega)}.
\end{eqnarray}
Integrating equation (\ref{eq:mm5}) from $0$ to $t$ gives
\begin{eqnarray}
\label{eq:mm6}
&&\|u^h\|^2_{L^2(\Omega)}+C\|\nabla u^h\|^p_{L^p(Q)}\nonumber\\&&\leq C\int_{0}^{t}\|u^h(x,s)\|^2_{L^2(\Omega)}ds+C\|u_0^h\|^2_{L^2(\Omega)}+C\|f\|^2_{L^2\big(0,T;L^2(\Omega)\big)}.
\end{eqnarray}
Ignoring the second term on the left side, since it is nonnegative, and applying Gronwall's lemma to (\ref{eq:mm6}), we obtain
\begin{eqnarray*}
\|u^h\|^2_{L^2(\Omega)}\leq C\|u_0^h\|^2_{L^2(\Omega)}+C\|f\|^2_{L^2\big(0,T;L^2(\Omega)\big)}.
\end{eqnarray*}
Using this estimate in (\ref{eq:mm6}) completes the proof of (\ref{reg_u}).
Finally, from equation (\ref{eq:set}), we obtain
\begin{eqnarray*}
\|y^h\|^2_{L^2(\Omega)}\leq C\|u_0^h\|^2_{L^2(\Omega)}+C\|f\|^2_{L^2\big(0,T;L^2(\Omega)\big)},
\end{eqnarray*}
where
\begin{equation*}
C=C\big(T,\|g\|_{L^\infty(0,T)},\|g'\|_{L^\infty(0,T)}\big),
\end{equation*}
as required.
\end{proof}

\begin{theorem}[Convergence]
Let $(u,y)$ and $(u^h,y^h)$ be solutions to problems (\ref{eq:ff27}) and (%
\ref{eq:ff28}), respectively. If $g,g^{\prime} \in L^{\infty}(0,T)$, $u_0\in
H^{r+1}(\Omega)$ and $f\in L^2(Q)$, then, for every $t\in [0,T]$,
\begin{eqnarray}  \label{eq:ff40}
\|u-u^h\|_{L^2(\Omega)}\leq Ch^{\frac{rp}{2(p-1)}}
\end{eqnarray}
and
\begin{eqnarray}  \label{eq:ff41}
\|y-y^h\|_{L^2(\Omega)}\leq Ch^{\frac{rp}{2(p-1)}},
\end{eqnarray}
where the constant $C$ does not depend on $h$ but may depend on $g$, $u$, $y$
and their derivatives.
\end{theorem}

\begin{proof}
Noting that $(u,y)$ and $(u^h,y^h)$ are solutions to problems (\ref{eq:ff27}) and (\ref{eq:ff28}), respectively, and subtracting the second equation of (\ref{eq:ff28}) from the second equation of (\ref{eq:ff27}), with $v = v^h\in\mathcal{S}^h\subset H_0^1$, we get
\begin{eqnarray*}
&&\int_{\Omega}^{}(y-y^h)v^hdx\nonumber\\&&=-\int_{\Omega}^{}v^h\int_{0}^{t}g(t-s)\big(y(x,s)-y^h(x,s)\big)dsdx+\int_{\Omega}^{}g(0)(u-u^h)v^hdx\nonumber\\&&-\int_{\Omega}^{}g(t)(u_0-u_0^h)v^hdx+\int_{\Omega}^{}v^h\int_{0}^{t}g'(t-s)\big(u(x,s)-u^h(x,s)\big)dsdx.
\end{eqnarray*}
By writing $(y-y^h)=(y-\Pi_hy)+(\Pi_hy-y^h)=\varphi+\psi$ and $(u-u^h)=(u-\Pi_hu)+(\Pi_hu-u^h)=\rho+\theta$, with $\Pi_hu$ the interpolation of $u$ in $\mathcal{S}^h$ , we have
\begin{eqnarray}
\label{eq:ff38}
&&\int_{\Omega}^{}\psi v^hdx\nonumber\\&&=-\int_{\Omega}^{}\varphi v^hdx-\int_{\Omega}^{}v^h\int_{0}^{t}g(t-s)\varphi(x,s)dsdx+g(0)\int_{\Omega}^{}\theta v^hdx\nonumber\\&&+g(0)\int_{\Omega}^{}\rho v^hdx-g(t)\int_{\Omega}^{}(u_0-u_0^h)v^hdx-\int_{\Omega}^{}v^h\int_{0}^{t}g(t-s)\psi(x,s)dsdx\nonumber\\&&+\int_{\Omega}^{}v^h\int_{0}^{t}g'(t-s)\rho(x,s)dsdx+\int_{\Omega}^{}v^h\int_{0}^{t}g'(t-s)\theta(x,s)dsdx.
\end{eqnarray}
Making $v^h=\psi\in\mathcal{S}^h$ and applying Young's inequality to (\ref{eq:ff38}), we have
\begin{eqnarray*}
&&\int_{\Omega}^{}\psi^2dx\nonumber\\&&\leq C(\epsilon)\int_{\Omega}^{}\varphi^2dx+\epsilon\int_{\Omega}^{}\psi^2dx+C(\epsilon)\int_{\Omega}^{}\theta^2dx+\epsilon\int_{\Omega}^{}\psi^2dx+ C(\epsilon)\int_{\Omega}^{}\rho^2dx\nonumber\\&&+\epsilon\int_{\Omega}^{}\psi^2dx+C(\epsilon)\int_{\Omega}^{}(u_0-u_0^h)^2dx+\epsilon\int_{\Omega}^{}\psi^2dx+\epsilon\int_{\Omega}^{}\psi^2dx\nonumber\\&&+C(\epsilon)\int_{\Omega}^{}\int_{0}^{t}\varphi^2dsdx+\epsilon\int_{\Omega}^{}\psi^2dx+C(\epsilon)\int_{\Omega}^{}\int_{0}^{t}\rho^2dsdx+\epsilon\int_{\Omega}^{}\psi^2dx\nonumber\\&&+C(\epsilon)\int_{\Omega}^{}\int_{0}^{t}\theta^2dsdx+\epsilon\int_{\Omega}^{}\psi^2dx+C(\epsilon)\int_{\Omega}^{}\int_{0}^{t}\psi^2(x,s)dsdx.
\end{eqnarray*}
Choosing $\epsilon$ appropriately and using Lemma \ref{eq:le3}, we obtain the inequality
\begin{eqnarray}
\label{eq:ff29}
\int_{\Omega}^{}\psi^2dx&\leq &Ch^{2(r+1)}+C\int_{\Omega}^{}\theta^2dx+C\int_{0}^{t}\int_{\Omega}^{}\theta^2(x,s)dxds\nonumber\\&&+C\int_{0}^{t}\int_{\Omega}^{}\psi^2(x,s)dxds.
\end{eqnarray}
Applying Gronwall's lemma to (\ref{eq:ff29}), we have
\begin{eqnarray}
\label{eq:set1}
\int_{\Omega}^{}\psi^2dx\leq Ch^{2(r+1)}+C\int_{\Omega}^{}\theta^2dx+C(1+t)\int_{0}^{t}\int_{\Omega}^{}\theta^2(x,s)dxds.
\end{eqnarray}
Now, subtracting the first equation in (\ref{eq:ff28}) from the first equation in (\ref{eq:ff27}), with $w=w^h\in\mathcal{S}^h$, we get
\begin{eqnarray*}
\label{eq:ff35}
&&\int_{\Omega}^{}\theta_t w^hdx+\int_{\Omega}^{}\left(\left|\Pi_h\nabla u\right|^{p-2}\Pi_h \nabla u-\left|\nabla u^h\right|^{p-2}\nabla u^h\right)\nabla w^hdx=\int_{\Omega}^{}\varphi w^hdx\nonumber\\&&+\int_{\Omega}^{}\psi w^hdx-\int_{\Omega}^{}\rho_tw^hdx+\int_{\Omega}^{}\left(\left|\Pi_h\nabla u\right|^{p-2}\Pi_h\nabla u-\left|\nabla u\right|^{p-2}\nabla u\right)\nabla w^hdx.
\end{eqnarray*}
Making $w^h=\theta\in\mathcal{S}^h$, applying Young's inequality and Lemma \ref{eq:le1}, we obtain
\begin{eqnarray*}
&&\frac{1}{2}\frac{d}{dt}\int_{\Omega}^{}\theta^2dx+C\int_{\Omega}^{}\left|\nabla\theta\right|^{p}dx\nonumber\\ &&\leq\frac{1}{2}\int_{\Omega}^{}\varphi^2dx+\frac{1}{2}\int_{\Omega}^{}\psi^2dx+\frac{1}{2}\int_{\Omega}^{}\theta^2dx+\frac{1}{2}\int_{\Omega}^{}\rho_t^2dx+\frac{1}{2}\int_{\Omega}^{}\theta^2dx\nonumber\\ &&+\epsilon\int_{\Omega}^{}|\nabla\theta|^pdx+C(\epsilon)\int_{\Omega}^{}|\nabla\rho|^{\frac{p}{p-1}}dx.
\end{eqnarray*}
Choosing $\epsilon$ appropriately and using estimate (\ref{eq:ff29}) and Lemma \ref{eq:le3}, we have
\begin{eqnarray}
\label{eq:ff39}
\frac{d}{dt}\int_{\Omega}^{}\theta^2dx\leq Ch^{2(r+1)}+Ch^{\frac{rp}{p-1}}+C\int_{\Omega}^{}\theta^2dx+C\int_{0}^{t}\int_{\Omega}^{}\theta^2(x,s)dxds.
\end{eqnarray}
Integrating (\ref{eq:ff39}) with respect to $t$, we get
\begin{eqnarray}
\label{eq:ff33}
\int_{\Omega}^{}\theta^2dx\leq Ch^{2(r+1)}+Ch^{\frac{rp}{p-1}}+C(1+t)\int_{0}^{t}\int_{\Omega}^{}\theta^2(x,s)dxds.
\end{eqnarray}
Applying Gronwall's lemma to (\ref{eq:ff33}), we obtain the inequality
\begin{eqnarray}
\label{eq:ff34}
\int_{\Omega}^{}\theta^2dx\leq Ch^{2(r+1)}+Ch^{\frac{rp}{p-1}}.
\end{eqnarray}
Equation (\ref{eq:ff34}) with the estimate for $\rho$ given in Lemma \ref{eq:le3} proves (\ref{eq:ff40}). Substituting equation (\ref{eq:ff34}) in equation (\ref{eq:set1}) and adding the estimate for $\varphi$, we obtain (\ref{eq:ff41}), as required.
\end{proof}

\section{Discretization in time}

The memory
term will be discretized using a numerical quadrature and, in order to maintain a good
convergence order, we use the Crank-Nicolson method along with the trapezoidal
quadrature.

\subsection{Crank-Nicolson method}

Consider the partition $0=t_0<t_1<\cdots <t_N=T$, with step $\delta=\frac{T}{%
N}$, of $[0,T]$. Evaluating (\ref{eq:ff28}), at $t={t_{k+\frac{1}{2}}}=\frac{%
t_{k+1}+t_k}{2}$, we obtain
\begin{eqnarray*}  \label{eq:ff42a}
&&\int_{\Omega}^{}u_t^h\Big(x,t_{k+\frac{1}{2}}\Big)w^hdx+\int_{\Omega}^{}%
\left|\nabla u^h\Big(x,t_{k+\frac{1}{2}}\Big)\right|^{p-2}\nabla u^{h}\Big(%
x,t_{k+\frac{1}{2}}\Big)\nabla w^hdx  \notag \\
&&=\int_{\Omega}^{}y^h\Big(x,t_{k+\frac{1}{2}}\Big)w^hdx+\int_{\Omega}^{}f%
\Big(x,t_{k+\frac{1}{2}}\Big)w^hdx
\end{eqnarray*}
and
\begin{eqnarray*}  \label{eq:ff43}
&&\int_{\Omega}^{}y^h\Big(x,t_{k+\frac{1}{2}}\Big)v^h(x)dx=-\int_{0}^{t_{k+%
\frac{1}{2}}}g\Big(t_{k+\frac{1}{2}}-s\Big)\int_{\Omega}^{}y^h(x,s)v^h(x)dxds
\notag \\
&&+g(0)\int_{\Omega}^{}u^h\Big(x,t_{k+\frac{1}{2}}\Big)v^h(x)dx-g\Big(t_{k+%
\frac{1}{2}}\Big)\int_{\Omega}^{}u_0(x)v^h(x)dx  \notag \\
&&+\int_{0}^{t_{k+\frac{1}{2}}}g^{\prime }\Big(t_{k+\frac{1}{2}}-s\Big)%
\int_{\Omega}^{}u^h(x,s)v^h(x)dxds  \notag \\
&&-\int_{0}^{t_{k+\frac{1}{2}}}g\Big(t_{k+\frac{1}{2}}-s\Big)%
\int_{\Omega}^{}f(x,s)v^h(x)dxds.
\end{eqnarray*}
Let us consider the approximations
\begin{equation*}
u_t^h\Big(x,t_{k+\frac{1}{2}}\Big)\approx \frac{u^h(x,t_{k+1})-u^h(x,t_{k})}{%
\delta}=\bar{\partial}u^{\big(k+\frac{1}{2}\big)},
\end{equation*}
\begin{equation*}
u^h\Big(x,t_{k+\frac{1}{2}}\Big)\approx\frac{u^h(x,t_{k+1})+u^h(x,t_{k})}{2}=%
\bar{u}^{\big(k+\frac{1}{2}\big)}
\end{equation*}
and
\begin{equation*}
y^h\Big(x,t_{k+\frac{1}{2}}\Big)\approx\frac{y^h(x,t_{k+1})+y^h(x,t_{k})}{2}=%
\bar{y}^{\big(k+\frac{1}{2}\big)}.
\end{equation*}
We will now approximate the
integrals in time using the trapezoidal quadrature so that the order of precision is maintained.
\begin{eqnarray*}
&&\int_{0}^{t_{k+\frac{1}{2}}}g^{\prime }\Big(t_{k+\frac{1}{2}}-s\Big)%
\int_{\Omega}^{}u^h(x,s)v^h(x)dxds  \notag \\
&&\approx\frac{\delta}{2}g^{\prime }\big(t_{k+\frac{1}{2}}\big)%
\int_{\Omega}^{}u_0^h(x)v^h(x)dx+\delta\sum\limits_{j=1}^{k-1}g^{\prime }%
\big(t_{k+\frac{1}{2}}-t_j\big)\int_{\Omega}^{}u^h(x,t_j)v^h(x)dx  \notag \\
&&+\frac{3\delta}{4}g^{\prime }\big(t_{k+\frac{1}{2}}-t_k\big)%
\int_{\Omega}^{}u^h(x,t_k)v^h(x)dx+\frac{\delta}{8}g^{\prime
}(0)\int_{\Omega}^{}u^h(x,t_k)v^h(x)dx  \notag \\
&&+\frac{\delta}{8}g^{\prime
}(0)\int_{\Omega}^{}u^h(x,t_{k+1})v^h(x)dx=Q_{g^{\prime }}(u^h)
\end{eqnarray*}
and
\begin{eqnarray*}
&&\int_{0}^{t_{k+\frac{1}{2}}}g\Big(t_{k+\frac{1}{2}}-s\Big)%
\int_{\Omega}^{}y^h(x,s)v^h(x)dxds  \notag \\
&& \approx\frac{\delta}{2}g\big(t_{k+\frac{1}{2}}\big)\int_{%
\Omega}^{}y_0^h(x)v^h(x)dx+\delta\sum\limits_{j=1}^{k-1}g\big(t_{k+\frac{1}{2%
}}-t_j\big)\int_{\Omega}^{}y^h(x,t_j)v^h(x)dx  \notag \\
&&+\frac{3\delta}{4}g\big(t_{k+\frac{1}{2}}-t_k\big)\int_{%
\Omega}^{}y^h(x,t_k)v^h(x)dx+\frac{\delta}{8}g(0)\int_{%
\Omega}^{}y^h(x,t_k)v^h(x)dx  \notag \\
&&+\frac{\delta}{8}g(0)\int_{\Omega}^{}y^h(x,t_{k+1})v^h(x)dx=Q_{g}(y^h).
\end{eqnarray*}

\subsection{Totally discrete formulation}

To simplify the notation, whenever there is no danger of confusion, we
will consider a function with the superscript $(j)$ to represent this function
evaluated at instant $t=t_j$.\newline
The totally discrete problem is to find $\big(U^{(k+1)},Y^{(k+1)}\big)$,
the solution to
\begin{eqnarray}  \label{eq:ff31}
\int_{\Omega}^{}\bar{\partial}U^{\big(k+\frac{1}{2}\big)}w^hdx
+\int_{\Omega}^{}\left|\nabla \bar{U}^{\big(k+\frac{1}{2}\big)%
}\right|^{p-2}\nabla \bar{U}^{\big(k+\frac{1}{2}\big)}\nabla
w^hdx=\int_{\Omega}^{}\bar{Y}^{\big(k+\frac{1}{2}\big)}w^hdx&&  \notag \\
+\int_{\Omega}^{}f^{\big(k+\frac{1}{2}\big)}w^hdx&&
\end{eqnarray}
and
\begin{eqnarray}  \label{eq:ff42}
\int_{\Omega}^{}\bar{Y}^{\big(k+\frac{1}{2}\big)}v^hdx =g(0)\int_{\Omega}^{}%
\bar{U}^{\big(k+\frac{1}{2}\big)}v^hdx -g\big(t_{k+\frac{1}{2}}\big)%
\int_{\Omega}^{}u_0^h(x)v^hdx&&  \notag \\
-Q_g(Y) +Q_{g^{\prime }}(U) -I(f),&&
\end{eqnarray}
where
\begin{equation*}  \label{eq:ff311}
I(f)=\int_{0}^{t_{k+\frac{1}{2}}}g\big(t_{k+\frac{1}{2}}-s\big)%
\int_{\Omega}^{}f(x,s)v^h(x)dxds.
\end{equation*}
Equation (\ref{eq:ff42}) is a linear system for $Y^{(k+1)}$. Obtaining a
solution $Y^{(k+1)}$ depending on $U^{(k+1)}$ from (\ref{eq:ff42}) and
substituting in (\ref{eq:ff31}),  one obtains a non linear equation for $%
U^{(k+1)}$. Then, using the fixed point theorem, it  follows easily that (\ref{eq:ff31}%
) has a solution.

\begin{theorem}[Uniqueness]
If $g,g^{\prime} \in L^{\infty}(0,T)$, the solution to the discrete problem is unique.
\end{theorem}

\begin{proof}
The proof is similar to that of Theorem \ref{teo_uni}, but it is more technical.
Suppose that $\big(U_1^{(k+1)},Y_1^{(k+1)}\big)$ and $\big(U_2^{(k+1)},Y_2^{(k+ 1)}\big)$ are two solutions to problem (\ref{eq:ff31})-(\ref{eq:ff42}). Subtracting the equation for $Y_2^{(k+1)}$ from the equation for $Y_1^{(k+1)}$, we obtain
\begin{eqnarray}
\label{eq:j15}
&&\left(\frac{1}{2}-\frac{\delta}{8}g(0)\right)\int_{\Omega}^{}\big(Y_1^{(k+1)}-Y_2^{(k+1)}\big)v^hdx\nonumber\\ &&=\left(\frac{3\delta}{4}g\big(t_{k+\frac{1}{2}}-t_k\big)+\frac{\delta}{8}g(0)-\frac{1}{2}\right)\int_{\Omega}^{}\big(Y_1^{(k)}-Y_2^{(k)}\big)v^hdx\nonumber\\
&&+\delta\sum\limits_{j=1}^{k-1}g\big(t_{k+\frac{1}{2}}-t_j\big)\int_{\Omega}^{}\big(Y_1^{(j)}-Y_2^{(j)}\big)v^hdx+\frac{\delta}{2}g\big(t_{k+\frac{1}{2}}\big)\int_{\Omega}^{}\big(Y_1^{(0)}-Y_2^{(0)}\big)v^hdx\nonumber\\
&&+\left(-\frac{1}{2}g(0)-\frac{\delta}{8}g'(0)\right)\int_{\Omega}^{}\big(U_1^{(k+1)}-U_2^{(k+1)}\big)v^hdx\nonumber\\
&&+\left(\frac{1}{2}g(0)-\frac{3\delta}{8}g'\big(t_{k+\frac{1}{2}}-t_k\big)-\frac{\delta}{8}g'(0)\right)\int_{\Omega}^{}\big(U_1^{(k)}-U_2^{(k)}\big)v^hdx\nonumber\\
&&+\delta\sum\limits_{j=1}^{k-1}g'\big(t_{k+\frac{1}{2}}-t_j\big)\int_{\Omega}^{}\big(U_1^{(j)}-U_2^{(j)}\big)v^hdx\nonumber\\
&&+\left(-g\big(t_{k+\frac{1}{2}}\big)+\frac{\delta}{8}g'\big(t_{k+\frac{1}{2}}\big)\right)\int_{\Omega}^{}\big(U_1^{(0)}-U_2^{(0)}\big)v^hdx.
\end{eqnarray}
Let $v^h=Y_1^{(k+1)}-Y_2^{(k+1)}$. As $g$ and $g'$ are bounded, we may apply Young's inequality in (\ref{eq:j15}), so
\begin{eqnarray}
\label{eq:j16}
&&\int_{\Omega}^{}\big(Y_1^{(k+1)}-Y_2^{(k+1)}\big)^2dx\leq C\int_{\Omega}^{}\big(Y_1^{(k)}-Y_2^{(k)}\big)^2dx\nonumber\\&& +C\delta\sum\limits_{j=1}^{k-1}\int_{\Omega}^{}\big(Y_1^{(j)}-Y_2^{(j)}\big)^2dx +C\int_{\Omega}^{}\big(U_1^{(k+1)}-U_2^{(k+1)}\big)^2dx\nonumber\\&& +C\int_{\Omega}^{}\big(U_1^{(k)}-U_2^{(k)}\big)^2dx +C\delta\sum\limits_{j=1}^{k-1}\int_{\Omega}^{}\big(U_1^{(j)}-U_2^{(j)}\big)^2dx.
\end{eqnarray}
Now, applying Gronwall's lemma to (\ref{eq:j16}), we have
\begin{eqnarray}
\label{eq:j17}
&\int_{\Omega}^{}\big(Y_1^{(k+1)}-Y_2^{(k+1)}\big)^2dx\leq& C\int_{\Omega}^{}\big(U_1^{(k+1)}-U_2^{(k+1)}\big)^2dx +C\int_{\Omega}^{}\big(U_1^{(k)}-U_2^{(k)}\big)^2dx\nonumber\\
&&+C\delta\sum\limits_{j=1}^{k-1}\int_{\Omega}^{}\big(U_1^{(j)}-U_2^{(j)}\big)^2dx.
\end{eqnarray}
Considering $w^h=\bar{U}^{\big(k+\frac{1}{2}\big)}$ in equation (\ref{eq:ff31}) and repeating the procedure for $ U_1^{(k+1)}$ and $U_2^{(k+1)}$, we obtain
\begin{eqnarray}
\label{eq:j18}
&&\int_{\Omega}^{}\big(U_1^{(k+1)}-U_2^{(k+1)}\big)^2dx -\int_{\Omega}^{}\big(U_1^{(k)}-U_2^{(k)}\big)^2dx\nonumber\\
&&+2\delta\int_{\Omega}^{}\left|\frac{\Big(\nabla U_1^{(k+1)}-\nabla U_2^{(k+1)}\Big)+\Big(\nabla U_1^{(k)}-\nabla U_2^{(k)}\Big)}{2}\right|^{p}dx\nonumber\\ &&=2\delta\int_{\Omega}^{}\Big(\bar{Y}_1^{\big(k+\frac{1}{2}\big)}-\bar{Y}_2^{\big(k+\frac{1}{2}\big)}\Big)\Big(\bar{U}_1^{\big(k+\frac{1}{2}\big)}-\bar{U}_2^{\big(k+\frac{1}{2}\big)}\Big)dx.
\end{eqnarray}
Applying Young's inequality to (\ref{eq:j18}) and (\ref{eq:j17}), we have
\begin{eqnarray*}
\label{eq:j20}
&\int_{\Omega}^{}\big(U_1^{(k+1)}-U_2^{(k+1)}\big)^2dx\leq& C\int_{\Omega}^{}\big(U_1^{(k)}-U_2^{(k)}\big)^2dx +C\int_{\Omega}^{}\big(U_1^{(k-1)}-U_2^{(k-1)}\big)^2dx\nonumber\\
&&+C\delta\sum_{j=1}^{k-1}\int_{\Omega}^{}\big(U_1^{(j)}-U_2^{(j)}\big)^2dx.
\end{eqnarray*}
As $\zeta^{(k+1)}=\int_{\Omega}^{}\big(U_1^{(k+1)}-U_2^{(k+1)}\big)^2dx\geq 0$ and $\zeta^{(0)}=0$, we may apply the discrete version of Gronwall's lemma. Therefore $\zeta^{(k)}=0$ for every $k\geq 0$ and so
\begin{eqnarray}
\label{eq:j21}
\int_{\Omega}^{}\big(U_1^{(k+1)}-U_2^{(k+1)}\big)^2dx=0,\quad U_1^{(k+1)}=U_2^{(k+1)}\;\,\mathrm{in}\;\,L^2(\Omega).
\end{eqnarray}
Going back to $Y$ and substituting (\ref{eq:j21}) in (\ref{eq:j17}), we have
\begin{equation*}
\int_{\Omega}^{}\big(Y_1^{(k+1)}-Y_2^{(k+1)}\big)^2dx=0,
\end{equation*}
which proves the required result.
\end{proof}

\begin{theorem}[Stability]
Let $\big(U^{(k+1)},Y^{(k+1)}\big)\in\mathcal{S}^h\times\mathcal{S}^h$ be
solutions to problem (\ref{eq:ff31}), (\ref{eq:ff42}). If $g,g^{\prime} \in L^{\infty}(0,T)$, $u_0\in L^2(\Omega)$ and $f\in L^2(Q)$, then, for every $k\geq
0$,
\begin{eqnarray*}
\|U^{(k+1)}\|_{L^2(\Omega)}^2\leq C\|u_0\|_{L^2(\Omega)}^2+C\|f\|_{L^2\big(%
0,T;L^2(\Omega)\big)}^2+C\delta\sum\limits_{j=0}^{k}\|f^{\big(j+\frac{1}{2}%
\big)}\|_{L^2(\Omega)}^2
\end{eqnarray*}
and
\begin{eqnarray*}
\|Y^{(k+1)}\|_{L^2(\Omega)}^2\leq C\|u_0\|_{L^2(\Omega)}^2+C\|f\|_{L^2\big(%
0,T;L^2(\Omega)\big)}^2+C\delta\sum\limits_{j=0}^{k}\|f^{\big(j+\frac{1}{2}%
\big)}\|_{L^2(\Omega)}^2,
\end{eqnarray*}
where
\begin{equation*}
C=C\big(T,\|g\|_{L^\infty(0,T)},\|g^{\prime }\|_{L^\infty(0,T)}\big).
\end{equation*}
\end{theorem}

\begin{proof}
We write equation (\ref{eq:ff42}) in the form
\begin{eqnarray}
\label{eq:j8}
&&\left(\frac{1}{2}-\frac{\delta}{8}g(0)\right)\int_{\Omega}^{}Y^{(k+1)}v^hdx\nonumber\\
&&=\left(\frac{3\delta}{4}g\big(t_{k+\frac{1}{2}}-t_k\big)+\frac{\delta}{8}g(0)-\frac{1}{2}\right)\int_{\Omega}^{}Y^{(k)}v^hdx\nonumber\\
&&+\delta\sum\limits_{j=1}^{k-1}g\big(t_{k+\frac{1}{2}}-t_j\big)\int_{\Omega}^{}Y^{(j)}v^hdx+\left(-\frac{1}{2}g(0)-\frac{\delta}{8}g'(0)\right)\int_{\Omega}^{}U^{(k+1)}v^hdx\nonumber\\
&&+\left(\frac{1}{2}g(0)-\frac{3\delta}{8}g'\big(t_{k+\frac{1}{2}}-t_k\big)-\frac{\delta}{8}g'(0)\right)\int_{\Omega}^{}U^{(k)}v^hdx\nonumber\\
&&+\delta\sum\limits_{j=1}^{k-1}g'\big(t_{k+\frac{1}{2}}-t_j\big)\int_{\Omega}^{}U^{(j)}v^hdx+\left(-g\big(t_{k+\frac{1}{2}}\big)+\frac{\delta}{8}g'\big(t_{k+\frac{1}{2}}\big)\right)\int_{\Omega}^{}u_{0}^hv^hdx\nonumber\\
&&-\int_{0}^{t_{k+\frac{1}{2}}}g\big(t_{k+\frac{1}{2}}-s\big)\int_{\Omega}^{}f(x,s)v^hdx,
\end{eqnarray}
and consider $v^h=Y^{(k+1)}$.  Then, since $g$ and $g'$ are bounded, we may apply Young's inequality to (\ref{eq:j8}) and thus obtain
\begin{eqnarray}
\label{eq:j9}
\|Y^{(k+1)}\|^2_{L^2(\Omega)}\leq C\|Y^{(k)}\|^2_{L^2(\Omega)}+C\delta\sum\limits_{j=1}^{k-1}\|Y^{(j)}\|^2_{L^2(\Omega)}+C\|U^{(k+1)}\|^2_{L^2(\Omega)}&&\nonumber\\
+C\|U^{(k)}\|^2_{L^2(\Omega)}+C\delta\sum\limits_{j=1}^{k-1}\|U^{(j)}\|^2_{L^2(\Omega)}+C\|u_0^h\|^2_{L^2(\Omega)}+C\|f\|^2_{L^2(\Omega)}.&&
\end{eqnarray}
Applying discrete version of Gronwall's lemma to (\ref{eq:j9}), we have
\begin{eqnarray}
\label{eq:j10}
&\|Y^{(k+1)}\|^2_{L^2(\Omega)}\leq &C\|U^{(k+1)}\|^2_{L^2(\Omega)}+C\|U^{(k)}\|^2_{L^2(\Omega)}+C\delta\sum\limits_{j=1}^{k-1}\|U^{(j)}\|^2_{L^2(\Omega)}\nonumber\\
&&+C\|f\|_{L^2\big(0,T;L^2(\Omega)\big)}^2+C\|u_0^h\|^2_{L^2(\Omega)}.
\end{eqnarray}
Returning now to equation (\ref{eq:ff31}) and considering  $w^h=\bar{U}^{\big(k+\frac{1}{2}\big)}$, we have
\begin{eqnarray}
\label{eq:j11}
&&\int_{\Omega}^{}\big(U^{(k+1)}\big)^2dx-\int_{\Omega}^{}\big(U^{(k)}\big)^2dx+2C\delta\int_{\Omega}^{}\left|\frac{\nabla U^{(k+1)}+\nabla U^{(k)}}{2}\right|^{p}dx\nonumber\\
&&=2\delta\int_{\Omega}^{}\bar{Y}^{\big(k+\frac{1}{2}\big)}\bar{U}^{\big(k+\frac{1}{2}\big)}dx+2\delta\int_{\Omega}^{}f^{\big(k+\frac{1}{2}\big)}\bar{U}^{\big(k+\frac{1}{2}\big)}dx.
\end{eqnarray}
Applying Young's inequality to (\ref{eq:j11}) and (\ref{eq:j10}), we get
\begin{eqnarray}
\label{eq:j13}
&\|U^{(k+1)}\|^2_{L^2(\Omega)} \leq & C\|U^{(k)}\|^2_{L^2(\Omega)}+C\|U^{(k-1)}\|^2_{L^2(\Omega)}+C\delta\sum\limits_{j=1}^{k-1}\|U^{(j)}\|^2_{L^2(\Omega)}\nonumber\\
&&+C\|u_0^h\|^2_{L^2(\Omega)}+C\|f\|_{L^2\big(0,T;L^2(\Omega)\big)}^2+C\delta\|f\|^2_{L^2(\Omega)}.
\end{eqnarray}
Applying the discrete version of Gronwall's lemma to (\ref{eq:j13}), we have
\begin{equation*}
\|U^{(k+1)}\|_{L^2(\Omega)}^2\leq C\|u_0^h\|_{L^2(\Omega)}^2+C\|f\|_{L^2\big(0,T;L^2(\Omega)\big)}^2+C\delta\sum\limits_{j=0}^{k}\|f^{\big(j+\frac{1}{2}\big)}\|_{L^2(\Omega)}^2.
\end{equation*}
Then, from equation (\ref{eq:j9}), we obtain
\begin{equation*}
\|Y^{(k+1)}\|_{L^2(\Omega)}^2\leq C\|u_0^h\|_{L^2(\Omega)}^2+C\|f\|_{L^2\big(0,T;L^2(\Omega)\big)}^2+C\delta\sum\limits_{j=0}^{k}\|f^{\big(j+\frac{1}{2}\big)}\|_{L^2(\Omega)}^2,
\end{equation*}
where
\begin{equation*}
C=C\big(T,\|g\|_{L^\infty(0,T)},\|g'\|_{L^\infty(0,T)}\big),
\end{equation*}
which proves the theorem.
\end{proof}

\begin{theorem}[Convergence]
Let $(u,y)$ and $(U,Y)$ be solutions to problems (\ref{eq:ff27}) and (\ref%
{eq:ff31})-(\ref{eq:ff42}), respectively. If $g,g^{\prime} \in L^{\infty}(0,T)$
and $u_0\in H^{r+1}(\Omega)$, then, for $\delta$ sufficiently small,
\begin{equation*}
\|u(x,t_k)-U^{(k)}(x)\|_{L^2(\Omega)}\leq C\Big(h^{r+1}+\delta^2+h^{\frac{rp%
}{2p-2}}+\delta^{\frac{p}{p-1}}\Big)
\end{equation*}
and
\begin{equation*}
\|y(x,t_k)-Y^{(k)}(x)\|_{L^2(\Omega)}\leq C\Big(h^{r+1}+\delta^2+h^{\frac{rp%
}{2p-2}}+\delta^{\frac{p}{p-1}}\Big),
\end{equation*}
where the constant $C$ does not depend on $h$ or $\delta$ but may depend on $%
g$, $u$, $y$ and their derivatives.
\end{theorem}

\begin{proof}
Subtracting (\ref{eq:ff42}) from the second equation of (\ref{eq:ff27}) evaluated at $t=t_{k+\frac{1}{2}}$ with $v=v^h$, we have
\begin{eqnarray*}
&&\int_{\Omega}^{}\Big(y^{\big(k+\frac{1}{2}\big)}-\bar{Y}^{\big(k+\frac{1}{2}\big)}\Big)v^hdx\nonumber\\ &&=g(0)\int_{\Omega}^{}\Big(u^{\big(k+\frac{1}{2}\big)}-\bar{U}^{\big(k+\frac{1}{2}\big)}\Big)v^h dx+g\big(t_{k+\frac{1}{2}}\big)\int_{\Omega}^{}\left(u^h_0-u_0\right)v^h dx\nonumber\\ &&-\int_{0}^{t_{k+\frac{1}{2}}}g\big(t_{k+\frac{1}{2}}-s\big)\int_{\Omega}^{}y(x,s)v^h(x)dxds+Q_g(Y)\nonumber\\ &&+\int_{0}^{t_{k+\frac{1}{2}}}g'\big(t_{k+\frac{1}{2}}-s\big)\int_{\Omega}^{}u(x,s)v^h(x)dxds-Q_{g'}(U)\nonumber\\ &&=g(0)I_1+g\Big(t_{k+\frac{1}{2}}\Big)I_2+I_3-I_4.
\end{eqnarray*}
Considering $y(x,t_j)-Y^{(j)}=\varphi(x,t_j)+\psi(x,t_j)=\varphi^{(j)}+\psi^{(j)}$ as before, we have
\begin{eqnarray*}
\int_{\Omega}^{}\Big(y^{\big(k+\frac{1}{2}\big)}-\bar{Y}^{\big(k+\frac{1}{2}\big)}\Big)v^hdx =\int_{\Omega}^{}\Big(y^{\big(k+\frac{1}{2}\big)}-\bar{y}^{\big(k+\frac{1}{2}\big)}\Big)v^hdx\nonumber&&\\
+\int_{\Omega}^{}\bar{\varphi}^{\big(k+\frac{1}{2}\big)}v^hdx+\int_{\Omega}^{}\bar{\psi}^{\big(k+\frac{1}{2}\big)}v^hdx.&&
\end{eqnarray*}
Letting $u(x,t_j)-U^{(j)}=\rho(x,t_j)+\theta(x,t_j)=\rho^{(j)}+\theta^{(j)}$, we have
\begin{eqnarray*}
I_1&=&\int_{\Omega}^{}\Big(u^{\big(k+\frac{1}{2}\big)}-\bar{U}^{\big(k+\frac{1}{2}\big)}\Big)v^hdx\nonumber\\ &=&\int_{\Omega}^{}\left(u^{\big(k+\frac{1}{2}\big)})-\bar{u}^{\big(k+\frac{1}{2}\big)}\right)v^hdx +\int_{\Omega}^{}\bar{\rho}^{\big(k+\frac{1}{2}\big)}v^hdx+\int_{\Omega}^{}\bar{\theta}^{\big(k+\frac{1}{2}\big)}v^hdx.
\end{eqnarray*}
We can write
\begin{eqnarray*}
\label{eq:ff44}
I_3=\int_{0}^{t_{k+\frac{1}{2}}}g\big(t_{k+\frac{1}{2}}-s\big)\int_{\Omega}^{}y(x,s)v^h(x)dxds-Q_g(y)+Q_g(y)-Q_g(Y),
\end{eqnarray*}
where
\begin{eqnarray*}
Q_g(y)-Q_g(Y)=\frac{\delta}{2}g\big(t_{k+\frac{1}{2}}\big)\int_{\Omega}^{}\big(y_0-Y^{(0)}\big)v^hdx +\delta\sum\limits_{j=1}^{k-1}g\big(t_{k+\frac{1}{2}}-t_j\big)\int_{\Omega}^{}\varphi^{(j)}v^hdx&&\nonumber\\
+\delta\sum\limits_{j=1}^{k-1}g\big(t_{k+\frac{1}{2}}-t_j\big)\int_{\Omega}^{}\psi^{(j)}v^hdx +\left(\frac{3\delta}{4}g\big(t_{k+\frac{1}{2}}-t_k\big)+\frac{\delta}{8}g(0)\right)\int_{\Omega}^{}\varphi^{(k)}v^hdx&&\nonumber\\
+\frac{3\delta}{4}g\big(t_{k+\frac{1}{2}}-t_k\big)\int_{\Omega}^{}\psi^{(k)}v^hdx +\frac{\delta}{8}g(0)\int_{\Omega}^{}\varphi^{(k+1)}v^hdx+\frac{\delta}{4}g(0)\int_{\Omega}^{}\bar{\psi}^{\big(k+\frac{1}{2}\big)}v^hdx.&&
\end{eqnarray*}
Likewise,
\begin{eqnarray*}
I_4=\int_{0}^{t_{k+\frac{1}{2}}}g'\big({t_{k+\frac{1}{2}}}-s\big)\int_{\Omega}^{}u(x,s)v^h(x)dxds-Q_{g'}(u)+Q_{g'}(u)-Q_{g'}(U),
\end{eqnarray*}
where
\begin{eqnarray*}
&&Q_{g'}(u)-Q_{g'}(U)=\frac{\delta}{2}g'\big(t_{k+\frac{1}{2}}\big)\int_{\Omega}^{}\big(u_0(x)-U^{(0)}\big)v^hdx\nonumber\\ &&+\delta\sum\limits_{j=1}^{k-1}g'\big(t_{k+\frac{1}{2}}-t_j\big)\int_{\Omega}^{}\rho^{(j)}v^hdx +\delta\sum\limits_{j=1}^{k-1}g'\big(t_{k+\frac{1}{2}}-t_j\big)\int_{\Omega}^{}\theta^{(j)}v^hdx\nonumber\\ &&+\left(\frac{3\delta}{4}g'\big(t_{k+\frac{1}{2}}-t_k\big)+\frac{\delta}{8}g'(0)\right)\int_{\Omega}^{}\rho^{(k)}v^hdx +\frac{3\delta}{4}g'\big(t_{k+\frac{1}{2}}-t_k\big)\int_{\Omega}^{}\theta^{(k)}v^hdx\nonumber\\ &&+\frac{\delta}{8}g'(0)\int_{\Omega}^{}\rho^{(k+1)}v^hdx +\frac{\delta}{4}g'(0)\int_{\Omega}^{}\bar{\theta}^{\big(k+\frac{1}{2}\big)}v^hdx,
\end{eqnarray*}
and hence we can write
\begin{eqnarray*}
&&\int_{\Omega}^{}\bar{\psi}^{\big(k+\frac{1}{2}\big)}v^hdx =-\int_{\Omega}^{}\left(y^{\big(k+\frac{1}{2}\big)}-\bar y^{\big(k+\frac{1}{2}\big)}\right)v^hdx -\frac{1}{2}\int_{\Omega}^{}\varphi^{(k+1)}v^hdx\nonumber\\
&& -\frac{1}{2}\int_{\Omega}^{}\varphi^{(k)}v^hdx +g(0)\int_{\Omega}^{}\left(u^{\big(k+\frac{1}{2}\big)}-\bar u^{\big(k+\frac{1}{2}\big)}\right)v^hdx +\frac{g(0)}{2}\int_{\Omega}^{}\rho^{(k+1)}v^hdx\nonumber\\ &&+\frac{g(0)}{2}\int_{\Omega}^{}\rho^{(k)}v^hdx +\int_{\Omega}^{}\bar{\theta}^{\big(k+\frac{1}{2}\big)}v^hdx +g\big(t_{k+\frac{1}{2}}\big)\int_{\Omega}^{}\big(U^{(0)}-u_0\big)v^hdx\nonumber\\ &&-\int_{0}^{t_{k+\frac{1}{2}}}g\big({t_{k+\frac{1}{2}}}-s\big)\int_{\Omega}^{}y(x,s)v^hdxds+Q_g(y)\nonumber\\ &&-\frac{\delta}{2}g\big(t_{k+\frac{1}{2}}\big)\int_{\Omega}^{}\big(y_0-Y^{(0)}\big)v^hdx -\delta\sum\limits_{j=1}^{k-1}g\big(t_{k+\frac{1}{2}}-t_j\big)\int_{\Omega}^{}\varphi^{(j)}v^hdx\nonumber\\ &&-\delta\sum\limits_{j=1}^{k-1}g\big(t_{k+\frac{1}{2}}-t_j\big)\int_{\Omega}^{}\psi^{(j)}v^hdx -\frac{3\delta}{4}g\big(t_{k+\frac{1}{2}}-t_k\big)\int_{\Omega}^{}\varphi^{(k)}v^hdx\nonumber\\ &&-\frac{3\delta}{4}g\big(t_{k+\frac{1}{2}}-t_k\big)\int_{\Omega}^{}\psi^{(k)}v^hdx -\frac{\delta}{8}g(0)\int_{\Omega}^{}\varphi^{(k)}v^hdx +\frac{\delta}{8}g(0)\int_{\Omega}^{}\varphi^{(k+1)}v^hdx\nonumber\\
&&-\frac{\delta}{4}g(0)\int_{\Omega}^{}\bar{\psi}^{\big(k+\frac{1}{2}\big)}v^hdx -\int_{0}^{t_{k+\frac{1}{2}}}g'\big({t_{k+\frac{1}{2}}}-s\big)\int_{\Omega}^{}u(x,s)v^h(x)dxds-Q_{g'}(u)\nonumber\\ &&+\frac{\delta}{2}g'\big(t_{k+\frac{1}{2}}\big)\int_{\Omega}^{}\big(u_0-U^{(0)}\big)v^h(x)dx +\delta\sum\limits_{j=1}^{k-1}g'\big(t_{k+\frac{1}{2}}-t_j\big)\int_{\Omega}^{}\rho^{(j)}v^hdx\nonumber\\ &&+\delta\sum\limits_{j=1}^{k-1}g'\big(t_{k+\frac{1}{2}}-t_j\big)\int_{\Omega}^{}\theta^{(j)}v^hdx +\frac{3\delta}{4}g'\big(t_{k+\frac{1}{2}}-t_k\big)\int_{\Omega}^{}\rho^{(k)}v^hdx\nonumber\\ &&+\frac{3\delta}{4}g'\big(t_{k+\frac{1}{2}}-t_k\big)\int_{\Omega}^{}\theta^{(k)}v^hdx +\frac{\delta}{8}g'(0)\int_{\Omega}^{}\rho^{(k)}v^hdx\nonumber\\ &&+\frac{\delta}{8}g'(0)\int_{\Omega}^{}\rho^{(k+1)}v^hdx +\frac{\delta}{4}g'(0)\int_{\Omega}^{}\bar{\theta}^{\big(k+\frac{1}{2}\big)}v^hdx.
\end{eqnarray*}
Using the estimates from the error term of the trapezoidal rule with $v^h=\bar{\psi}^{\big(k+\frac{1}{2}\big)}$ and applying Young's inequality, we have
\begin{eqnarray}
&&\left(1-25\epsilon+\frac{\delta g(0)}{4}\right)\|\bar{\psi}^{\big(k+\frac{1}{2}\big)}\|^2_{L^2(\Omega)} \leq C(\epsilon)\|y^{\big(k+\frac{1}{2}\big)}-\bar y^{\big(k+\frac{1}{2}\big)}\|^2_{L^2(\Omega)}\nonumber\\ &&+C(\epsilon)(1+\delta^2)\|\varphi^{(k+1)}\|^2_{L^2(\Omega)} +C(\epsilon)(1+2\delta^2)\|\varphi^{(k)}\|^2_{L^2(\Omega)} \nonumber\\
&&+C(\epsilon)\|u^{\big(k+\frac{1}{2}\big)}-\bar u^{\big(k+\frac{1}{2}\big)}\|^2_{L^2(\Omega)} +C(\epsilon)(1+\delta^2)\|\rho^{(k+1)}\|^2_{L^2(\Omega)}
\nonumber\\
&&+C(\epsilon)(1+2\delta^2)\|\rho^{(k)}\|^2_{L^2(\Omega)} +C(\epsilon)(1+\delta^2)\|u_0(x)-U^{(0)}\|^2_{L^2(\Omega)}
\nonumber\\ &&+C(\epsilon)\delta^4\left(\|y\|^2_{L^\infty\big(0,T;L^2(\Omega)\big)} +\|y_t\|^2_{L^\infty\big(0,T;L^2(\Omega)\big)}
+\|y_{tt}\|^2_{L^\infty\big(0,T;L^2(\Omega)\big)}\right)\nonumber\\ &&+C(\epsilon)\delta^2\|y_0-Y^{(0)}\|^2_{L^2(\Omega)} +C(\epsilon)\delta^2\sum\limits_{j=1}^{k-1}\|\theta^{(j)}\|^2_{L^2(\Omega)} \nonumber\\ &&+C(\epsilon)\delta^2\sum\limits_{j=1}^{k-1}\|\psi^{(j)}\|^2_{L^2(\Omega)} +C(\epsilon)\delta^2\sum\limits_{j=1}^{k-1}\|\varphi^{(j)}\|^2_{L^2(\Omega)} +C(\epsilon)\delta^2\|\psi^{(k)}\|^2_{L^2(\Omega)}\nonumber\\ &&+C(\epsilon)\delta^4\left(\|u\|^2_{L^\infty\big(0,T;L^2(\Omega)\big)} +\|u_t\|^2_{L^\infty\big(0,T;L^2(\Omega)\big)}
+\|u_{tt}\|^2_{L^\infty\big(0,T;L^2(\Omega)\big)}\right)\nonumber\\ &&+C(\epsilon)\delta^2\sum\limits_{j=1}^{k-1}\|\rho^{(j)}\|^2_{L^2(\Omega)} +C(\epsilon)\delta^2\|\theta^{(k)}\|^2_{L^2(\Omega)} \nonumber\\ &&+C(\epsilon)(1+\delta^2)\|\bar{\theta}^{\big(k+\frac{1}{2}\big)}\|^2_{L^2(\Omega)}. \label{eq_g0}
\end{eqnarray}
By Lemma \ref{eq:le3}, some interpolation error bounds and choosing appropriate $\epsilon$ and $\delta$, we can write
\begin{eqnarray*}
\label{eq:ff45}
\|\bar{\psi}^{\big(k+\frac{1}{2}\big)}\|^2_{L^2(\Omega)}\leq C\big(h^{2(r+1)}+\delta^4\big)+C\|\bar{\theta}^{\big(k+\frac{1}{2}\big)}\|^2_{L^2(\Omega)} +C\delta^2\sum\limits_{j=0}^{k}\|\theta^{(j)}\|^2_{L^2(\Omega)}&&\nonumber\\ +C\delta^2\sum\limits_{j=0}^{k}\|\psi^{(j)}\|^2_{L^2(\Omega)}.
\end{eqnarray*}
Then
\begin{eqnarray*}
\label{eq:ff46}
\|\psi^{(k+1)}\|^2_{L^2(\Omega)}\leq  C\|\psi^{(k)}\|^2_{L^2(\Omega)}+C\big(h^{2(r+1)}+\delta^4\big)+C\|\bar{\theta}^{\big(k+\frac{1}{2}\big)}\|^2_{L^2(\Omega)}&&\nonumber\\
+C\delta^2\sum\limits_{j=0}^{k}\|\theta^{(j)}\|^2_{L^2(\Omega)}+C\delta^2\sum\limits_{j=0}^{k}\|\psi^{(j)}\|^2_{L^2(\Omega)}.&&
\end{eqnarray*}
By the discrete version of Gronwall's lemma,
\begin{eqnarray*}
\|\psi^{(k+1)}\|^2_{L^2(\Omega)}\leq C(1+\delta)\big(h^{2(r+1)}+\delta^4\big)
+C\|\theta^{(k+1)}\|^2_{L^2(\Omega)} +C\|\theta^{(k)}\|^2_{L^2(\Omega)}&&\nonumber\\
+C\|\theta^{(k-1)}\|^2_{L^2(\Omega)} +C(1+\delta)\delta^2\sum\limits_{j=0}^{k-1}\|\theta^{(j)}\|^2_{L^2(\Omega)}.&&
\end{eqnarray*}
Now we go back to the equation for $u$. Subtracting equation (\ref{eq:ff31}) from the first equation of (\ref{eq:ff27}), evaluated at $t=t_{k+\frac{1}{2}}$, and considering $w=w ^h\in \mathcal{S}^h$, we obtain
\begin{eqnarray*}
&&\int_{\Omega}^{}\left(u_t^{\big(k+\frac{1}{2}\big)}-\bar{\partial}U^{\big(k+\frac{1}{2}\big)}\right)w^hdx\nonumber\\
&&+\int_{\Omega}^{}\left(\left|\nabla u^{\big(k+\frac{1}{2}\big)}\right|^{p-2}\nabla u^{\big(k+\frac{1}{2}\big)}-\left|\nabla \bar{U}^{\big(k+\frac{1}{2}\big)}\right|^{p-2}\nabla \bar{U}^{\big(k+\frac{1}{2}\big)}\right)\nabla w^hdx\nonumber\\ &&=\int_{\Omega}^{}\left(y^{\big(k+\frac{1}{2}\big)}-\bar{Y}^{\big(k+\frac{1}{2}\big)}\right)w^hdx.
\end{eqnarray*}
In this case,
\begin{eqnarray*}
&&\int_{\Omega}^{}\left(u_t^{\big(k+\frac{1}{2}\big)}-\bar{\partial}u^{\big(k+\frac{1}{2}\big)}\right)w^hdx +\int_{\Omega}^{}\left(\bar{\partial}u^{\big(k+\frac{1}{2}\big)} -\bar{\partial}U^{\big(k+\frac{1}{2}\big)}\right)w^hdx\nonumber\\
&&+\int_{\Omega}^{}\left(\left|\nabla u^{\big(k+\frac{1}{2}\big)}\right|^{p-2}\nabla u^{\big(k+\frac{1}{2}\big)}-\left|\nabla\overline{\Pi_h u}^{\big(k+\frac{1}{2}\big)}\right|^{p-2}\nabla\overline{\Pi_h u}^{\big(k+\frac{1}{2}\big)}\right) \nabla w^h dx\nonumber\\
&&+\int_{\Omega}^{}\left(\left|\nabla\overline{\Pi_h u}^{\big(k+\frac{1}{2}\big)}\right|^{p-2}\nabla\overline{\Pi_h u}^{\big(k+\frac{1}{2}\big)}-\left|\nabla \bar{U}^{\big(k+\frac{1}{2}\big)}\right|^{p-2}\nabla \bar{U}^{\big(k+\frac{1}{2}\big)}\right)\nabla w^h dx\nonumber\\
&&=\int_{\Omega}^{}\left(y^{\big(k+\frac{1}{2}\big)}-\bar y^{\big(k+\frac{1}{2}\big)}\right)w^hdx +\int_{\Omega}^{}\left(\bar y^{\big(k+\frac{1}{2}\big)}-\bar{Y}^{\big(k+\frac{1}{2}\big)}\right)w^hdx,
\end{eqnarray*}
and therefore
\begin{eqnarray*}
&&\int_{\Omega}^{}\bar{\partial}\theta^{\big(k+\frac{1}{2}\big)}w^hdx\nonumber\\
&&+\int_{\Omega}^{}\left(\left|\nabla\overline{\Pi_h u}^{\big(k+\frac{1}{2}\big)}\right|^{p-2}\nabla\overline{\Pi_h u}^{\big(k+\frac{1}{2}\big)}-\left|\nabla \bar{U}^{\big(k+\frac{1}{2}\big)}\right|^{p-2}\nabla \bar{U}^{\big(k+\frac{1}{2}\big)}\right)\nabla w^hdx\nonumber\\
&&=-\int_{\Omega}^{}\left(u_t^{\big(k+\frac{1}{2}\big)} -\bar{\partial}u^{\big(k+\frac{1}{2}\big)}\right)w^hdx -\int_{\Omega}^{}\bar{\partial}\rho^{\big(k+\frac{1}{2}\big)}w^hdx\nonumber\\
&&+\int_{\Omega}^{}\left(\left|\nabla u^{\big(k+\frac{1}{2}\big)}\right|^{p-2}\nabla u^{\big(k+\frac{1}{2}\big)}-\left|\nabla\overline{\Pi_h u}^{\big(k+\frac{1}{2}\big)}\right|^{p-2}\nabla\overline{\Pi_h u}^{\big(k+\frac{1}{2}\big)}\right) \nabla w^h dx\nonumber\\ &&+\int_{\Omega}^{}\left(y^{\big(k+\frac{1}{2}\big)}-\bar y^{\big(k+\frac{1}{2}\big)}\right)w^hdx +\int_{\Omega}^{}\bar{\varphi}^{\big(k+\frac{1}{2}\big)}w^hdx +\int_{\Omega}^{}\bar{\psi}^{\big(k+\frac{1}{2}\big)}w^hdx.
\end{eqnarray*}
Taking $w^h=\bar{\theta}^{\big(k+\frac{1}{2}\big)}$ and applying Young's inequality, we obtain
\begin{eqnarray*}
&&\int_{\Omega}^{}\bar{\partial}\theta^{\big(k+\frac{1}{2}\big)}\bar{\theta}^{\big(k+\frac{1}{2}\big)}dx+C_2\|\nabla \bar{\theta}^{\big(k+\frac{1}{2}\big)}\|^p_{L^p(\Omega)}\nonumber\\
&&\leq C_{\epsilon_1}\|u_t^{\big(k+\frac{1}{2}\big)}-\bar{\partial}u^{\big(k+\frac{1}{2}\big)}\|^2_{L^2(\Omega)} +\epsilon_1\|\bar{\theta}^{\big(k+\frac{1}{2}\big)}\|^2_{L^2(\Omega)} +C_{\epsilon_1}\|\bar{\partial}\rho^{\big(k+\frac{1}{2}\big)}\|^2_{L^2(\Omega)}\nonumber\\
&& +\epsilon_1\|\bar{\theta}^{\big(k+\frac{1}{2}\big)}\|^2_{L^2(\Omega)}+\epsilon_2\|\nabla \bar{\theta}^{\big(k+\frac{1}{2}\big)}\|^p_{L^p(\Omega)}\nonumber\\
&&+C_{\epsilon_2}\left\|\left|\nabla u^{\big(k+\frac{1}{2}\big)}\right|^{p-2}\nabla u^{\big(k+\frac{1}{2}\big)}-\left|\nabla \bar u^{\big(k+\frac{1}{2}\big)}\right|^{p-2}\nabla \bar u^{\big(k+\frac{1}{2}\big)}\right\|^{\frac{p}{p-1}}_{L^2(\Omega)}\nonumber\\
&&+C_{\epsilon_2}\left\|\left|\nabla \bar u^{\big(k+\frac{1}{2}\big)}\right|^{p-2}\nabla \bar u^{\big(k+\frac{1}{2}\big)} -\left|\nabla\overline{\Pi_h u}^{\big(k+\frac{1}{2}\big)}\right|^{p-2}\nabla\overline{\Pi_h u}^{\big(k+\frac{1}{2}\big)}\right\|^{\frac{p}{p-1}}_{L^2(\Omega)}\nonumber\\
&&+\epsilon_2\|\nabla \bar{\theta}^{\big(k+\frac{1}{2}\big)}\|^p_{L^2(\Omega)} +\epsilon_1\|\bar{\theta}^{\big(k+\frac{1}{2}\big)}\|^2_{L^2(\Omega)} +C_{\epsilon_1}\|y^{\big(k+\frac{1}{2}\big)}-\bar y^{\big(k+\frac{1}{2}\big)}\|^2_{L^2(\Omega)}\nonumber\\
&&+C_{\epsilon_1}\|\bar{\varphi}^{\big(k+\frac{1}{2}\big)}\|^{2}_{L^2(\Omega)} +\epsilon_1\|\bar{\theta}^{\big(k+\frac{1}{2}\big)}\|^2_{L^2(\Omega)}+C_{\epsilon_1}\|\bar{\psi}^{\big(k+\frac{1}{2}\big)}\|^2_{L^2(\Omega)}+\epsilon_1\|\bar{\theta}^{\big(k+\frac{1}{2}\big)}\|^2_{L^2(\Omega)}.
\end{eqnarray*}
Using Lemma \ref{eq:le1} and Lemma \ref{eq:le3}, some numerical differentiation and interpolation error bounds and choosing $\epsilon_2<\frac{C_2}{2}$, we can write
\begin{eqnarray*}
&&\frac{1}{2\delta}\left(\|\theta^{(k+1)}\|^2_{L^2(\Omega)}+\|\theta^{(k)}\|^2_{L^2(\Omega)}\right)\nonumber\\&&\leq C\delta^4+C\|\theta^{(k+1)}\|^2_{L^2(\Omega)}+C\|\theta^{(k)}\|^2_{L^2(\Omega)}+Ch^{2(r+1)}+C\delta^{\frac{2p}{p-1}}\nonumber\\ &&+Ch^{\frac{rp}{p-1}}+C\delta^4+Ch^{2(r+1)}+C\|\psi^{(k+1)}\|^2_{L^2(\Omega)}+C\|\psi^{(k)}\|^2_{L^2(\Omega)},
\end{eqnarray*}
that is,
\begin{eqnarray*}
(1-C\delta)\|\theta^{(k+1)}\|^2_{L^2(\Omega)}\leq C\delta\big(h^{2(r+1)}+\delta^4\big)+C\delta\Big(\delta^{\frac{2p}{p-1}}+h^{\frac{rp}{p-1}}\Big)&&\nonumber\\
+(C\delta-1)\|\theta^{(k)}\|^2_{L^2(\Omega)}+C\delta\|\theta^{(k-1)}\|^2_{L^2(\Omega)}+C\delta^3\sum\limits_{j=0}^{k-1}\|\theta^{(j)}\|^2_{L^2(\Omega)}.&&
\end{eqnarray*}
For $\delta$ sufficiently small, we have
\begin{eqnarray*}
\|\theta^{(k+1)}\|^2_{L^2(\Omega)}\leq C\delta\big(h^{2(r+1)}+\delta^4\big)+C\delta\Big(\delta^{\frac{2p}{p-1}}+h^{\frac{rp}{p-1}}\Big)-C\|\theta^{(k)}\|^2_{L^2(\Omega)}&&\nonumber\\ +C\delta\|\theta^{(k-1)}\|^2_{L^2(\Omega)}+C\delta^3\sum\limits_{j=0}^{k-1}\|\theta^{(j)}\|^2_{L^2(\Omega)}.&&
\end{eqnarray*}
From the discrete version of Gronwall's lemma,
\begin{eqnarray*}
\|\theta^{(k+1)}\|^2_{L^2(\Omega)}\leq C\delta\left(h^{2(r+1)}+\delta^4+h^{\frac{rp}{p-1}}+\delta^{\frac{2p}{p-1}}\right).
\end{eqnarray*}
Returning to the equation of $\|\psi^{(k+1)}\|^2_{L^2(\Omega)}$,
\begin{eqnarray*}
\|\psi^{(k+1)}\|^2_{L^2(\Omega)}\leq C\big(h^{2(r+1)}+\delta^4\big)+C\delta\left(h^{2(r+1)}+\delta^4+h^{\frac{rp}{p-1}}+\delta^{\frac{2p}{p-1}}\right)&&\nonumber\\
+C(1+\delta)\delta^2\left(h^{2(r+1)}+\delta^4+h^{\frac{rp}{p-1}}+\delta^{\frac{2p}{p-1}}\right).&&
\end{eqnarray*}
Finally, adding the estimates of $\rho^{(k+1)}$ and $\varphi^{(k+1)}$ given by Lemma \ref{eq:le3}, the required result is obtained.
\end{proof}

We notice that in (\ref{eq_g0}) if $g(0)\geq 0$ then $\delta$ can be any positive value, otherwise should be $\delta<\frac{-4}{g(0)}$ sufficiently small, for example $\delta=\frac{-1}{g(0)}$ for $\epsilon=\frac{1}{100}$.

\section{Final comments}

In this paper, we applied the finite element method with a polynomial basis of
degree $r$ complemented with the Crank-Nicolson method and the trapezoid
quadrature to a class of evolution differential equations with p-Laplacian
and memory. The memory term was separated from the $p$-Laplacian using a
mixed formulation. We demonstrated the existence, uniqueness and regularity
of the discrete solutions under mild conditions on the data. We also
obtained the convergence order depending on $p$ in the classical norms. It
was found that the convergence order decreases as $p\to\infty$, but it is
always bigger than $\frac{r}{2}$ for $h$ and bigger than 1 for $\delta$.

As future work, we intend to find an efficient method to solve the nonlinear system of algebraic
equations and to implement the method in a
computational system, such as in a Matlab environment, to illustrate the theoretical results. An interesting challenge is to do a similar study for equation (%
\ref{eq:30}) with $p$ depending on $x$. The fact that the $p(x)$-Laplacian
is not homogeneous makes the problem more complicated than the problem with
constant $p$.

\section*{Funding}

This work was partially supported by FEDER through the - Programa
Operacional Factores de Competitividade, FCT - Funda\c{c}\~{a}o para a Ci%
\^{e}ncia e a Tecnologia [Grant N. UIDB/00212/2020] and Santander [Grant N.
BID/ICI-FC/Santander Universidades-UBI/2015].



\begin{thebibliography}{26}
\expandafter\ifx\csname natexlab\endcsname\relax\def\natexlab#1{#1}\fi
\providecommand{\bibinfo}[2]{#2}
\ifx\xfnm\relax \def\xfnm[#1]{\unskip,\space#1}\fi
\bibitem[{Antontsev and Shmarev(2015)}]{AS_book}
\bibinfo{author}{S.~Antontsev}, \bibinfo{author}{S.~Shmarev},
  \bibinfo{title}{Evolution {PDE}s with nonstandard growth conditions},
  volume~\bibinfo{volume}{4} of \textit{\bibinfo{series}{Atlantis Studies in
  Differential Equations}}, \bibinfo{publisher}{Atlantis Press, Paris},
  \bibinfo{year}{2015}. \bibinfo{note}{Existence, uniqueness, localization,
  blow-up}.
\bibitem[{Antontsev et~al.(2016)Antontsev, Shmarev, Simsen and
  Simsen}]{MR3462616}
\bibinfo{author}{S.~Antontsev}, \bibinfo{author}{S.~Shmarev},
  \bibinfo{author}{J.~Simsen}, \bibinfo{author}{M.S. Simsen},
  \bibinfo{title}{On the evolution {$p$}-{L}aplacian with nonlocal memory},
  \bibinfo{journal}{Nonlinear Anal.} \bibinfo{volume}{134}
  (\bibinfo{year}{2016}) \bibinfo{pages}{31--54}.
\bibitem[{Antontsev et~al.(2019)Antontsev, Shmarev, Simsen and
  Stefanello~Simsen}]{MR3919666}
\bibinfo{author}{S.~Antontsev}, \bibinfo{author}{S.~Shmarev},
  \bibinfo{author}{J.~Simsen}, \bibinfo{author}{M.~Stefanello~Simsen},
  \bibinfo{title}{Differential inclusion for the evolution {$p(x)$}-{L}aplacian
  with memory}, \bibinfo{journal}{Electron. J. Differential Equations}
  (\bibinfo{year}{2019}) \bibinfo{pages}{Paper No. 26, 28}.
\bibitem[{Barbu and Malik(1979)}]{BM79}
\bibinfo{author}{V.~Barbu}, \bibinfo{author}{M.A. Malik},
  \bibinfo{title}{Semilinear integro-differential equations in {H}ilbert
  space}, \bibinfo{journal}{J. Math. Anal. Appl.} \bibinfo{volume}{67}
  (\bibinfo{year}{1979}) \bibinfo{pages}{452--475}.
\bibitem[{Barrett and Liu(1994)}]{BL94}
\bibinfo{author}{J.W. Barrett}, \bibinfo{author}{W.B. Liu},
  \bibinfo{title}{Finite element approximation of the parabolic
  {$p$}-{L}aplacian}, \bibinfo{journal}{SIAM J. Numer. Anal.}
  \bibinfo{volume}{31} (\bibinfo{year}{1994}) \bibinfo{pages}{413--428}.
\bibitem[{Chen and Shih(1998)}]{chen_book}
\bibinfo{author}{C.~Chen}, \bibinfo{author}{T.~Shih}, \bibinfo{title}{Finite
  element methods for integrodifferential equations},
  volume~\bibinfo{volume}{9} of \textit{\bibinfo{series}{Series on Applied
  Mathematics}}, \bibinfo{publisher}{World Scientific Publishing Co., Inc.,
  River Edge, NJ}, \bibinfo{year}{1998}.
\bibitem[{Chipot(2000)}]{MR1801735}
\bibinfo{author}{M.~Chipot}, \bibinfo{title}{Elements of nonlinear analysis},
  Birkh\"{a}user Advanced Texts: Basler Lehrb\"{u}cher. [Birkh\"{a}user
  Advanced Texts: Basel Textbooks], \bibinfo{publisher}{Birkh\"{a}user Verlag,
  Basel}, \bibinfo{year}{2000}.
\bibitem[{Chipot and Savitska(2014)}]{MR3250760}
\bibinfo{author}{M.~Chipot}, \bibinfo{author}{T.~Savitska},
  \bibinfo{title}{Nonlocal {$p$}-{L}aplace equations depending on the {$L^p$}
  norm of the gradient}, \bibinfo{journal}{Adv. Differential Equations}
  \bibinfo{volume}{19} (\bibinfo{year}{2014}) \bibinfo{pages}{997--1020}.
\bibitem[{Chow(1989)}]{Cho89}
\bibinfo{author}{S.S. Chow}, \bibinfo{title}{Finite element error estimates for
  nonlinear elliptic equations of monotone type}, \bibinfo{journal}{Numer.
  Math.} \bibinfo{volume}{54} (\bibinfo{year}{1989}) \bibinfo{pages}{373--393}.
\bibitem[{Ciarlet(2002)}]{MR1930132}
\bibinfo{author}{P.G. Ciarlet}, \bibinfo{title}{The finite element method for
  elliptic problems}, volume~\bibinfo{volume}{40} of
  \textit{\bibinfo{series}{Classics in Applied Mathematics}},
  \bibinfo{publisher}{Society for Industrial and Applied Mathematics (SIAM),
  Philadelphia, PA}, \bibinfo{year}{2002}. \bibinfo{note}{Reprint of the 1978
  original [North-Holland, Amsterdam; MR0520174 (58 \#25001)]}.
\bibitem[{Crandall et~al.(1978)Crandall, Londen and Nohel}]{CLN78}
\bibinfo{author}{M.G. Crandall}, \bibinfo{author}{S.O. Londen},
  \bibinfo{author}{J.A. Nohel}, \bibinfo{title}{An abstract nonlinear
  {V}olterra integrodifferential equation}, \bibinfo{journal}{J. Math. Anal.
  Appl.} \bibinfo{volume}{64} (\bibinfo{year}{1978}) \bibinfo{pages}{701--735}.
\bibitem[{DiBenedetto(1993)}]{DiB_book}
\bibinfo{author}{E.~DiBenedetto}, \bibinfo{title}{Degenerate parabolic
  equations}, Universitext, \bibinfo{publisher}{Springer-Verlag, New York},
  \bibinfo{year}{1993}.
\bibitem[{Evans(1998)}]{MR1625845}
\bibinfo{author}{L.C. Evans}, \bibinfo{title}{Partial differential equations},
  volume~\bibinfo{volume}{19} of \textit{\bibinfo{series}{Graduate Studies in
  Mathematics}}, \bibinfo{publisher}{American Mathematical Society, Providence,
  RI}, \bibinfo{year}{1998}.
\bibitem[{Glowinski and Marrocco(1975)}]{GM75}
\bibinfo{author}{R.~Glowinski}, \bibinfo{author}{A.~Marrocco},
  \bibinfo{title}{Sur l'approximation, par \'{e}l\'{e}ments finis d'ordre un,
  et la r\'{e}solution, par p\'{e}nalisation-dualit\'{e}, d'une classe de
  probl\`emes de {D}irichlet non lin\'{e}aires}, \bibinfo{journal}{Rev.
  Fran\c{c}aise Automat. Informat. Recherche Op\'{e}rationnelle S\'{e}r. Rouge
  Anal. Num\'{e}r.} \bibinfo{volume}{9} (\bibinfo{year}{1975})
  \bibinfo{pages}{41--76}.
\bibitem[{Himadan(2021)}]{MR4220412}
\bibinfo{author}{A.~Himadan}, \bibinfo{title}{Well defined extinction time of
  solutions for a class of weak-viscoelastic parabolic equation with positive
  initial energy}, \bibinfo{journal}{AIMS Math.} \bibinfo{volume}{6}
  (\bibinfo{year}{2021}) \bibinfo{pages}{4331--4344}.
\bibitem[{Lions(1969)}]{Lio_book}
\bibinfo{author}{J.L. Lions}, \bibinfo{title}{Quelques m\'{e}thodes de
  r\'{e}solution des probl\`emes aux limites non lin\'{e}aires},
  \bibinfo{publisher}{Dunod; Gauthier-Villars, Paris}, \bibinfo{year}{1969}.
\bibitem[{MacCamy(1976)}]{Cam76}
\bibinfo{author}{R.C. MacCamy}, \bibinfo{title}{Stability theorems for a class
  of functional differential equations}, \bibinfo{journal}{SIAM J. Appl. Math.}
  \bibinfo{volume}{30} (\bibinfo{year}{1976}) \bibinfo{pages}{557--576}.
\bibitem[{Mustapha et~al.(2011)Mustapha, Brunner, Mustapha and
  Sch\"{o}tzau}]{MBMS11}
\bibinfo{author}{K.~Mustapha}, \bibinfo{author}{H.~Brunner},
  \bibinfo{author}{H.~Mustapha}, \bibinfo{author}{D.~Sch\"{o}tzau},
  \bibinfo{title}{An {$hp$}-version discontinuous {G}alerkin method for
  integro-differential equations of parabolic type}, \bibinfo{journal}{SIAM J.
  Numer. Anal.} \bibinfo{volume}{49} (\bibinfo{year}{2011})
  \bibinfo{pages}{1369--1396}.
\bibitem[{Nohel(1981)}]{Noh79}
\bibinfo{author}{J.A. Nohel}, \bibinfo{title}{Nonlinear {V}olterra equations
  for heat flow in materials with memory}, in: \bibinfo{booktitle}{Integral and
  functional differential equations ({P}roc. {C}onf., {W}est {V}irginia
  {U}niv., {M}organtown, {W}. {V}a., 1979)}, volume~\bibinfo{volume}{67} of
  \textit{\bibinfo{series}{Lecture Notes in Pure and Appl. Math.}},
  \bibinfo{publisher}{Dekker, New York}, \bibinfo{year}{1981}, pp.
  \bibinfo{pages}{3--82}.
\bibitem[{Pani et~al.(0708)Pani, Fairweather and Fernandes}]{PFF08}
\bibinfo{author}{A.K. Pani}, \bibinfo{author}{G.~Fairweather},
  \bibinfo{author}{R.I. Fernandes}, \bibinfo{title}{Alternating direction
  implicit orthogonal spline collocation methods for an evolution equation with
  a positive-type memory term}, \bibinfo{journal}{SIAM J. Numer. Anal.}
  \bibinfo{volume}{46} (\bibinfo{year}{2007/08}) \bibinfo{pages}{344--364}.
\bibitem[{Reddy et~al.(2019)Reddy, Sinha and Cuminato}]{RSC19}
\bibinfo{author}{G.M.M. Reddy}, \bibinfo{author}{R.K. Sinha},
  \bibinfo{author}{J.A. Cuminato}, \bibinfo{title}{A posteriori error analysis
  of the {C}rank-{N}icolson finite element method for parabolic
  integro-differential equations}, \bibinfo{journal}{J. Sci. Comput.}
  \bibinfo{volume}{79} (\bibinfo{year}{2019}) \bibinfo{pages}{414--441}.
\bibitem[{Sinha et~al.(2006)Sinha, Ewing and Lazarov}]{SEL06}
\bibinfo{author}{R.K. Sinha}, \bibinfo{author}{R.E. Ewing},
  \bibinfo{author}{R.D. Lazarov}, \bibinfo{title}{Some new error estimates of a
  semidiscrete finite volume element method for a parabolic
  integro-differential equation with nonsmooth initial data},
  \bibinfo{journal}{SIAM J. Numer. Anal.} \bibinfo{volume}{43}
  (\bibinfo{year}{2006}) \bibinfo{pages}{2320--2343}.
\bibitem[{Sinha et~al.(2009)Sinha, Ewing and Lazarov}]{SEL09}
\bibinfo{author}{R.K. Sinha}, \bibinfo{author}{R.E. Ewing},
  \bibinfo{author}{R.D. Lazarov}, \bibinfo{title}{Mixed finite element
  approximations of parabolic integro-differential equations with nonsmooth
  initial data}, \bibinfo{journal}{SIAM J. Numer. Anal.} \bibinfo{volume}{47}
  (\bibinfo{year}{2009}) \bibinfo{pages}{3269--3292}.
\bibitem[{Tchier et~al.(2021)Tchier, Dassios, Tawfiq and Ragoub}]{TDTR21}
\bibinfo{author}{F.~Tchier}, \bibinfo{author}{I.~Dassios},
  \bibinfo{author}{F.~Tawfiq}, \bibinfo{author}{L.~Ragoub}, \bibinfo{title}{On
  the approximate solution of partial integro-differential equations using the
  pseudospectral method based on chebyshev cardinal functions},
  \bibinfo{journal}{Mathematics} \bibinfo{volume}{9} (\bibinfo{year}{2021})
  \bibinfo{pages}{286}.
\bibitem[{Wang and Hong(2019)}]{WH19}
\bibinfo{author}{W.~Wang}, \bibinfo{author}{Q.~Hong}, \bibinfo{title}{Two-grid
  economical algorithms for parabolic integro-differential equations with
  nonlinear memory}, \bibinfo{journal}{Appl. Numer. Math.}
  \bibinfo{volume}{142} (\bibinfo{year}{2019}) \bibinfo{pages}{28--46}.
\bibitem[{Zennir and Miyasita(2020)}]{ZM20}
\bibinfo{author}{K.~Zennir}, \bibinfo{author}{T.~Miyasita},
  \bibinfo{title}{Lifespan of solutions for a class of pseudo-parabolic
  equation with weak-memory}, \bibinfo{journal}{Alexandria Engineering Journal}
  \bibinfo{volume}{59} (\bibinfo{year}{2020}) \bibinfo{pages}{957--964}.

\end{thebibliography}
\end{document}